\documentclass[11pt]{article}
\usepackage{amscd}
\usepackage{amssymb,amsfonts,amsmath,psfrag,amscd,cite}
\usepackage{amsmath}
  \usepackage{paralist}
  \usepackage{graphics}
  \usepackage{epsfig}
  \usepackage[colorlinks=true]{hyperref}
  \usepackage[all]{xy}
\usepackage{graphicx}
\newtheorem{theo}{Theorem}

\newtheorem{rema}[theo]{Remark}

\makeatletter \@addtoreset{equation}{section}
\@addtoreset{theo}{section} \makeatother

\begin{document}
\date{}
\title{Dynamical Equations of Controlled Rigid Spacecraft with a Rotor}
\author{Hong Wang  \\
School of Mathematical Sciences and LPMC,\\
Nankai University, Tianjin 300071, P.R.China\\
E-mail: hongwang@nankai.edu.cn\\
May 4, 2020} \maketitle

{\bf Abstract.} In this paper, we consider the controlled rigid spacecraft with
an internal rotor as a regular point reducible regular controlled
Hamiltonian (RCH) system. In the cases of coincident and
non-coincident centers of buoyancy and gravity, we first give
the regular point reduction and the dynamical vector field
of the reduced controlled rigid spacecraft-rotor system, respectively.
Then, we derive precisely the geometric constraint
conditions of the reduced symplectic form for the
dynamical vector field of the regular point reducible controlled spacecraft-rotor system,
that is, the two types of Hamilton-Jacobi equation for the reduced
controlled spacecraft-rotor system by calculation in
detail. These researches reveal the deeply internal
relationships of the geometrical structures of phase spaces, the dynamical
vector fields and controls of the system. \\

{\bf Keywords:}\; controlled rigid spacecraft with a rotor, \;\;
regular controlled Hamiltonian system, \;\; coincident and non-coincident
centers, \;\; regular point reduction, \;\; Hamilton-Jacobi equation.\\

{\bf AMS Classification:} 70H33, \;\; 70H20, \;\; 70Q05.

\tableofcontents

\section{Introduction}

A regular controlled Hamiltonian (RCH) system is a
Hamiltonian system with external force and control,
which is defined in Marsden et al.\cite{mawazh10}. In general,
an RCH system, under the actions of external force and control, is not
Hamiltonian, however, it is a dynamical system closely related to a
Hamiltonian system, and it can be explored and studied by extending
the methods for external force and control in the study of Hamiltonian systems.
Thus, we can emphasize explicitly the impact of external force
and control in the study for the RCH systems. In particular,
in Marsden et al.\cite{mawazh10}, the authors
we give the regular point reduction and
the regular orbit reduction for an RCH system with
symmetry, by analyzing carefully the geometrical and topological
structures of the phase space and the reduced phase space of the
corresponding Hamiltonian system.These research work not only gave a
variety of reduction methods for the RCH systems, but
also showed a variety of relationships of controlled Hamiltonian
equivalence of these systems. Now, it is a natural problem
if there is a practical RCH system and how to show the effect on controls
in regular symplectic reductions of the system.
In this paper, we consider the rigid spacecraft-rotor system
with the control torque acting on an internal rotor
as a regular point reducible RCH system on the generalization
of a Lie group, and as the application of the above theoretical result,
we first give the regular point reduced controlled spacecraft-rotor
systems by calculation in detail, in the cases of coincident and
non-coincident centers of buoyancy and gravity.\\

We note also that Hamilton-Jacobi theory for
a Hamiltonian system is an important research subject.
On the one hand, it provides
a characterization of the generating functions of certain
time-dependent canonical transformations, such that a given
Hamiltonian system in such a form that its solutions are extremely
easy to find by reduction to the equilibrium. 
On the other hand, it is possible in many cases that
Hamilton-Jacobi equation provides an immediate way to integrate the
equation of motion of system, even when the problem of Hamiltonian
system itself has not been or cannot be solved completely, see Abraham and
Marsden \cite{abma78}, Arnold \cite{ar89} and Marsden and Ratiu
\cite{mara99}. In addition, the Hamilton-Jacobi equation is also 
fundamental in the study of the quantum-classical relationship in quantization, 
and it also plays an important role in the study
of stochastic dynamical systems, see Woodhouse \cite{wo92}, Ge and
Marsden \cite{gema88} and
L\'{a}zaro-Cam\'{i} and Ortega \cite{laor09}. For these reasons it
is described as a useful tools in the study of Hamiltonian system
theory, and has been extensively developed in past many years and
become one of the most active subjects in the study of modern
applied mathematics and analytical mechanics. \\

Let $Q$ be a smooth manifold and $TQ$ the tangent bundle, $T^* Q$
the cotangent bundle with the canonical symplectic form $\omega$,
and the projection $\pi_Q: T^* Q \rightarrow Q $ induces the map $
T\pi_{Q}: TT^* Q \rightarrow TQ. $
From Abraham and Marsden in \cite{abma78}, we know that the following
theorem give a classical description of Hamilton-Jacobi equation
from the generating function and the geometrical point of view.
\begin{theo}
Assume that the triple $(T^*Q,\omega,H)$ is a Hamiltonian system
with Hamiltonian vector field $X_H$, and $W: Q\rightarrow
\mathbb{R}$ is a given generating function. Then the following two assertions
are equivalent:\\
\noindent $(\mathrm{i})$ For every curve $\sigma: \mathbb{R}
\rightarrow Q $ satisfying $\dot{\sigma}(t)= T\pi_Q
(X_H(\mathbf{d}W(\sigma(t))))$, $\forall t\in \mathbb{R}$, then
$\mathbf{d}W \cdot \sigma $ is an integral curve of the Hamiltonian
vector field $X_H$.\\
\noindent $(\mathrm{ii})$ $W$ satisfies the Hamilton-Jacobi equation
$H(q^i,\frac{\partial W}{\partial q^i})=E, $ where $E$ is a
constant.
\end{theo}

From the above theorem, we know that
the assertion $(\mathrm{i})$ with equivalent to
Hamilton-Jacobi equation by the generating function
gives a geometric constraint condition of the canonical symplectic form
on the cotangent bundle $T^*Q$
for Hamiltonian vector field of the system.
Thus, the Hamilton-Jacobi equation reveals the deeply internal relationships of
the generating function, the canonical symplectic form
and the dynamical vector field of a Hamiltonian system.
But, from Marsden et al.\cite{mawazh10} we know that,
since the symplectic reduced system of a
Hamiltonian system with symmetry defined on the cotangent bundle
$T^*Q$ may not be a Hamiltonian system on a cotangent bundle,
then we cannot give the Hamilton-Jacobi theorem for the Marsden-Weinstein
reduced Hamiltonian system just like same as the above theorem.
We have to look for a new way. It is worthy of noting that,
in Wang \cite{wa17}, the two new types of Hamilton-Jacobi
equations for Hamiltonian system and regular reduced Hamiltonian
systems are given. by using the (reduced) symplectic
forms and the (reduced) dynamical vector fields, which are the
development of classical Hamilton-Jacobi equation given by Abraham and Marsden
\cite{abma78}. \\

Since an RCH system defined on the cotangent bundle
$T^*Q$, in general,  may not be a Hamiltonian system,
and it has yet no generating function,
we cannot give the Hamilton-Jacobi theorem for the RCH system
and its regular reduced systems just like same as the above theorem.
But, in Wang \cite{wa13d} the author can give precisely
the geometric constraint conditions
of the (reduced) symplectic forms for the
dynamical vector fields of an RCH system and its regular reduced systems,
that is, two types of Hamilton-Jacobi equations, which are the
development of the above two types of
Hamilton-Jacobi equations for a Hamiltonian system and its 
Marsden-Weinstein reduced Hamiltonian system given in Wang
\cite{wa17}. In this paper,
as the application of the above theoretical result,
we derive precisely the geometric constraint
conditions of the reduced symplectic forms for the
dynamical vector fields of the regular point reducible controlled spacecraft-rotor system,
that is, the two types of Hamilton-Jacobi equations for the reduced
controlled spacecraft-rotor system by calculation in
detail, in the cases of coincident and
non-coincident centers of buoyancy and gravity.\\

A brief of outline of this paper is as follows. In the second
section, we first review some relevant basic facts
about rigid spacecraft with an internal rotor, and give the
Hamiltonian function of the controlled rigid spacecraft-rotor system,
in the cases of coincident and non-coincident centers of buoyancy and gravity,
respectively, which will be used in subsequent sections.
As the application of the theoretical result of regular
point reduction of an RCH system given by Marsden et al
\cite{mawazh10}, in the third section we consider the controlled rigid
spacecraft-rotor system with the control torque acting on an internal rotor
as a regular point reducible RCH system on
the generalization of rotation group $\textmd{SO}(3)\times S^1$ and
on the generalization of Euclidean group $\textmd{SE}(3)\times S^1$,
respectively, we give the regular point reduced controlled spacecraft-rotor
systems, in the cases of coincident and non-coincident
centers of buoyancy and gravity. Moreover, as an application of the
Hamilton-Jacobi theoretical result for the regular reduced RCH
system given by Wang \cite{wa13d}, in the fourth section, we
 derive precisely the two types of Hamilton-Jacobi equations for
 the regular point reduced controlled rigid spacecraft-rotor
systems by calculation in detail, in the cases of coincident and non-coincident centers
of buoyancy and gravity. These research work reveal the deeply internal
relationships of the geometrical structures of phase spaces, the dynamical
vector fields and controls of the controlled rigid spacecraft-rotor system,
and develop the application of the regular symplectic reduction
and Hamilton-Jacobi theory for the RCH systems
with symmetries, and make us have much deeper understanding and
recognition for the structure of Hamiltonian systems and RCH
systems.

\section{The Rigid Spacecraft with an Internal Rotor}

In this section, we first give the Hamiltonian of rigid spacecraft with an internal rotor,
in the cases of coincident and non-coincident centers of buoyancy and gravity,
respectively. We first review some relevant 
basic facts about rigid spacecraft with an internal
rotor, which will be used in subsequent sections. We shall follow
the notations and conventions introduced in Marsden
\cite{ma92}, Marsden and Ratiu \cite{mara99}, and Marsden et al
\cite{mawazh10}. In this paper, we assume that all manifolds are
real, smooth and finite dimensional and all actions are smooth left
actions. For convenience, we also assume that all controls appearing
in this paper are the admissible controls.

\subsection{With Coincident Centers of Buoyancy and Gravity}

We first describe a rigid spacecraft carrying an internal "non-mass" rotor,
which is called a carrier body, where "non-mass" means that the mass of a rotor
is very very small comparing with the mass of the rigid spacecraft.
We first assume that the external forces
and torques acting on the rigid spacecraft-rotor system are due to
buoyancy and gravity. In general, it is possible that the
spacecraft's center of buoyancy may not be coincident with its
center of gravity. But, in this subsection we assume that the
spacecraft is symmetric and to have uniformly distributed mass, and
the center of buoyancy and the center of gravity are coincident.
Denote by $O$ the center of mass of the system in the carrier body frame and
at $O$ place a set of (orthogonal) body axes. Assume that the body
coordinate axes are aligned with principal axes of the carrier body,
and the rotor is aligned along the third principal axis, see Marsden
\cite{ma92}. Next, we assume that the rotor spins under the influence of a control torque $u$
acting on the rotor. If translations are ignored and only rotations
are considered, in this case, then the configuration space is
$Q=\textmd{SO}(3)\times S^1$, with the first factor being the
attitude of rigid spacecraft and the second factor being the angle
of rotor. The corresponding phase space is the cotangent bundle
$T^*Q$ and locally, $T^\ast Q=T^\ast \textmd{SO}(3)\times T^\ast S^1$, where
$T^\ast S^1 \cong T^\ast \mathbb{R}$ locally, with the
canonical symplectic form $\omega_Q$.
By using the local left trivialization, locally,
$T^\ast \textmd{SO}(3)\cong \textmd{SO}(3)\times \mathfrak{so}^\ast(3)$
and $T^*\mathbb{R} \cong \mathbb{R}\times \mathbb{R}^{*}$,
then we have that locally,
$T^*Q \cong \textmd{SO}(3)\times \mathfrak{so}^\ast(3)
\times \mathbb{R}\times \mathbb{R}^{*}$.
For convenience, in the following we denote uniformly that, locally,
$Q= \textmd{SO}(3)\times \mathbb{R}, $ and
$T^* Q= T^*(\textmd{SO}(3)\times \mathbb{R})\cong \textmd{SO}(3)\times \mathfrak{so}^\ast(3)
\times \mathbb{R} \times \mathbb{R}^{*}$.\\

Let $I=diag(I_1,I_2,I_3)$ be the matrix of inertia moment of the carrier
body in the body fixed frame, which is a principal body frame,
and $J_3$ be the moment of
inertia of rotor around its rotation axis. Let $J_{3k},\; k=1,2,$ be
the moments of inertia of the rotor around the $k$th principal axis
with $k=1,2,$ and denote by $\bar{I}_k= I_k +J_{3k}, \; k=1,2, \;
\bar{I}_3= I_3. $ Let $\Omega=(\Omega_1,\Omega_2,\Omega_3)$ be the
angular velocity vector of the rigid spacecraft-rotor system computed with respect to the axes
fixed in the carrier body and $(\Omega_1,\Omega_2,\Omega_3)\in
\mathfrak{so}(3)$. Let $\alpha$ be the relative angle of rotor and
$\dot{\alpha}$ the relative angular velocity vector of rotor about the third
principal axis with respect to a carrier body fixed frame.
For convenience, we assume the total mass of the system $m=1$.\\

Now, by the local left trivialization, locally,
$T\textmd{SO}(3)\cong \textmd{SO}(3)\times \mathfrak{so}(3)$
and $T\mathbb{R}\cong \mathbb{R}\times \mathbb{R}$,
then we have that locally,
$TQ \cong \textmd{SO}(3)\times \mathfrak{so}(3)
\times \mathbb{R}\times \mathbb{R}$.
We consider the Lagrangian of the rigid spacecraft-rotor system
$L(A,\Omega,\alpha,\dot{\alpha}):
TQ\cong \textmd{SO}(3)\times\mathfrak{so}(3)\times
\mathbb{R}\times\mathbb{R}\to \mathbb{R}$, which is the total
kinetic energy of the rigid spacecraft plus the kinetic energy of
rotor, given by
$$L(A,\Omega,\alpha,\dot{\alpha})=\dfrac{1}{2}[\bar{I}_1\Omega_1^2
+\bar{I}_2\Omega_2^2+\bar{I}_3\Omega_3^2
+J_3(\Omega_3+\dot{\alpha})^2],$$
where $A\in \textmd{SO}(3)$,
$\Omega=(\Omega_1,\Omega_2,\Omega_3)\in \mathfrak{so}(3)$, $\alpha
\in \mathbb{R}$, $\dot{\alpha} \in \mathbb{R}$. If we introduce the
conjugate angular momentum, given by
$$ \Pi_k= \dfrac{\partial
L}{\partial \Omega_k} =\bar{I}_k\Omega_k, \; k=1,2, \;\; \Pi_3=
\dfrac{\partial L}{\partial \Omega_3}
=\bar{I}_3\Omega_3+J_3(\Omega_3+\dot{\alpha}), \;\;
 l = \dfrac{\partial L}{\partial \dot{\alpha}}
=J_3(\Omega_3+\dot{\alpha}),$$
and by the Legendre transformation
\begin{align*}
FL: TQ \cong \textmd{SO}(3)\times\mathfrak{so}(3)\times\mathbb{R}\times\mathbb{R}& \to
T^*Q \cong \textmd{SO}(3)\times \mathfrak{so}^\ast(3)\times\mathbb{R}\times\mathbb{R}^*, \\
(A,\Omega,\alpha,\dot{\alpha}) &\to (A,\Pi, \alpha, l),
\end{align*}
where $\Pi=(\Pi_1,\Pi_2,\Pi_3)\in \mathfrak{so}^\ast(3)$, $l \in
\mathbb{R}^*$, we have the Hamiltonian $H(A,\Pi,\alpha,l): T^* Q \cong
\textmd{SO}(3)\times \mathfrak{so}^\ast
(3)\times\mathbb{R}\times\mathbb{R}^*\to \mathbb{R}$ given by
\begin{align}
H(A,\Pi,\alpha, l) &= \Omega\cdot \Pi+\dot{\alpha}\cdot
l-L(A,\Omega,\alpha,\dot{\alpha}) \nonumber \\
&= \frac{1}{2}[\frac{\Pi_1^2}{\bar{I}_1}+\frac{\Pi_2^2}{\bar{I}_2}
+\frac{(\Pi_3-l)^2}{\bar{I}_3}+\frac{l^2}{J_3}].
\end{align}
In this case, in order to give the dynamical vector field and the two types of
Hamilton-jacobi equation of the controlled rigid spacecraft-rotor
system, we need to consider the regular point reduction of the controlled spacecraft-rotor system and
give precisely the $R_p$-reduced symplectic form of the cotangent bundle $T^* Q\cong
\textmd{SO}(3)\times \mathfrak{so}^*(3)\times\mathbb{R}\times\mathbb{R}^{*}$.

\subsection{With Non-Coincident Centers of Buoyancy and Gravity}

Since it is possible that the rigid spacecraft's center of buoyancy may
not be coincident with its center of gravity, in this subsection
then we consider the rigid spacecraft-rotor system with non-coincident
centers of buoyancy and gravity. We fix an orthogonal coordinate
frame to the carrier body with origin located at the center of
buoyancy and axes aligned with the principal axes of the carrier
body, and the rotor is aligned along the third principal axis, see
Marsden \cite{ma92}, and Leonard and Marsden \cite{lema97}. Moreover, assume that when
the carrier body is oriented so that the carrier body frame is aligned
with the inertial frame, the third principal axis aligns with the
direction of gravity. The vector from the center of buoyancy to the
center of gravity with respect to the carrier body frame is $h\chi$,
where $\chi$ is an unit vector on the line connecting the two
centers which is assumed to be aligned along the third principal
axis, and $h$ is the length of this segment. Assume that the total mass of the carrier
body $m=1$, and the magnitude of gravitational acceleration
is denoted $g$, and let $\Gamma$ be the unit vector viewed by an observer
moving with the carrier body, and the rotor
spins under the influence of a control torque $u$ acting on the rotor..
In this case, the configuration space is
$Q=\textmd{SO}(3)\circledS \mathbb{R}^3\times S^1 \cong
\textmd{SE}(3)\times S^1$, with the first factor being the attitude
of the rigid spacecraft and the drift of the rigid spacecraft in the rotational
process and the second factor being the angle of rotor. The
corresponding phase space is the cotangent bundle $T^*Q$ and locally, $T^\ast Q =T^\ast
\textmd{SE}(3)\times T^\ast S^1$, where $T^\ast S^1 \cong T^\ast
\mathbb{R}$ locally, with the canonical symplectic form $\omega_Q$.
By using the local left trivialization, locally,
$T^\ast \textmd{SE}(3)\cong \textmd{SE}(3)\times \mathfrak{se}^\ast(3)$
and $T^*\mathbb{R}\cong \mathbb{R}\times \mathbb{R}^{*}$,
then we have that locally,
$T^*Q \cong \textmd{SE}(3)\times \mathfrak{se}^\ast(3)
\times \mathbb{R}\times \mathbb{R}^{*}$.
For convenience, in the following we denote uniformly that, locally,
$Q= \textmd{SE}(3)\times \mathbb{R}, $ and
$T^* Q= T^*(\textmd{SE}(3)\times \mathbb{R})\cong \textmd{SE}(3)\times \mathfrak{se}^\ast(3)
\times \mathbb{R} \times \mathbb{R}^{*}$.\\

Now, by the local left trivialization, locally,
$T\textmd{SE}(3)\cong \textmd{SE}(3)\times \mathfrak{se}(3)$
and $T\mathbb{R}\cong \mathbb{R}\times \mathbb{R}$,
then we have that locally,
$TQ \cong \textmd{SE}(3)\times \mathfrak{se}(3)
\times \mathbb{R}\times \mathbb{R}$.
We consider the Lagrangian of the rigid spacecraft-rotor system
$L(A,c,\Omega,\Gamma,\alpha,\dot{\alpha}): TQ \cong
\textmd{SE}(3)\times\mathfrak{se}(3)\times\mathbb{R}\times\mathbb{R}\to
\mathbb{R}$, which is the total kinetic energy of the rigid
spacecraft plus the kinetic energy of rotor minus potential energy
of the rigid spacecraft-rotor system, given by
$$L(A,c,\Omega,\Gamma,\alpha,\dot{\alpha})=\dfrac{1}{2}[\bar{I}_1\Omega_1^2
+\bar{I}_2\Omega_2^2+\bar{I}_3\Omega_3^2
+J_3(\Omega_3+\dot{\alpha})^2]- gh\Gamma \cdot \chi,$$ where
$(A,c)\in \textmd{SE}(3)$, $\Omega=(\Omega_1,\Omega_2,\Omega_3)\in
\mathfrak{so}(3)$, $\alpha \in \mathbb{R}$, $\dot{\alpha} \in
\mathbb{R}$, and the variable $\Gamma \in \mathbb{R}^3$ is regarded as a
parameter with respect to potential energy of the system,
$(\Omega, \Gamma)\in \mathfrak{se}(3)$.
If we introduce the conjugate angular momentum, given
by $$\Pi_k= \dfrac{\partial L}{\partial \Omega_k} =\bar{I}_k\Omega_k,
\; k=1,2, \;\; \Pi_3= \dfrac{\partial L}{\partial \Omega_3}
=\bar{I}_3\Omega_3+J_3(\Omega_3+\dot{\alpha}),
\;\; l =\dfrac{\partial L}{\partial \dot{\alpha}}
=J_3(\Omega_3+\dot{\alpha}),$$ and by the Legendre transformation
with the parameter $\Gamma$, that is,
\begin{align*}
FL: TQ \cong\textmd{SE}(3)\times\mathfrak{se}(3)\times\mathbb{R}\times\mathbb{R}
& \to T^*Q \cong\textmd{SE}(3)\times
\mathfrak{se}^\ast(3)\times\mathbb{R}\times\mathbb{R}^*, \\
(A,c,\Omega,\Gamma,\alpha,\dot{\alpha}) & \to (A,c,\Pi,\Gamma,\alpha, l),
\end{align*}
where $\Pi=(\Pi_1,\Pi_2,\Pi_3)\in \mathfrak{so}^\ast(3)$, $(\Pi,
\Gamma)\in \mathfrak{se}^\ast(3)$, $l \in \mathbb{R}^*$, we have the
Hamiltonian \\ $H(A,c,\Pi,\Gamma,\alpha,l): T^* Q \cong
\textmd{SE}(3)\times \mathfrak{se}^\ast
(3)\times\mathbb{R}\times\mathbb{R}^* \to \mathbb{R}$ given by
\begin{align}
 H(A,c,\Pi,\Gamma,\alpha,l) &=\Omega\cdot \Pi+\dot{\alpha}\cdot
l-L(A,c,\Omega,\Gamma,\alpha,\dot{\alpha}) \nonumber \\
&=\frac{1}{2}[\frac{\Pi_1^2}{\bar{I}_1}+\frac{\Pi_2^2}{\bar{I}_2}
+\frac{(\Pi_3-l)^2}{\bar{I}_3}+\frac{l^2}{J_3}]+ gh\Gamma \cdot
\chi \; . \end{align}
In this case, in order to give the dynamical vector field and the two types of
Hamilton-jacobi equation of the controlled rigid spacecraft-rotor
system, we need to consider the regular point reduction of the controlled spacecraft-rotor system and
give precisely the $R_p$-reduced symplectic form of the cotangent bundle $T^* Q\cong
\textmd{SE}(3)\times \mathfrak{se}^*(3)\times\mathbb{R}\times\mathbb{R}^{*}$.

\section{Symmetric Reduction of the Controlled Rigid Spacecraft-Rotor System}

We know that the main goal of
reduction theory in mechanics is to use conservation laws and the
associated symmetries to reduce the number of dimensions of a
mechanical system required to be described. So, such reduction
theory is regarded as a useful tool for simplifying and studying
concrete mechanical systems. In particular, the Marsden-Weinstein reduction
for a Hamiltonian system with symmetry and momentum map
is famous work, and great developments have been
obtained around the work in the theoretical study and applications
of mathematics, mechanics and physics. See
Abraham and Marsden \cite{abma78}, Abraham et al.
\cite{abmara88}, Arnold \cite{ar89}, Libermann and Marle
\cite{lima87}, Marsden \cite{ma92},  Marsden et al.
\cite{mamiorpera07, mamora90}, Marsden and Perlmutter \cite{mape00},
Marsden and Ratiu \cite{mara99},
Marsden and Weinstein \cite{mawe74}, Meyer \cite{me73},
Nijmeijer and Van der Schaft \cite {nivds90}
and Ortega and Ratiu \cite{orra04}.\\

It is worthy of noting that the authors in Marsden et al. \cite{mawazh10}
set up the regular reduction theory for the RCH systems
with symplectic structures and symmetries on a symplectic fiber
bundle, as an extension of the  Marsden-Weinstein reduction theory of
Hamiltonian systems under regular controlled Hamiltonian equivalence
conditions, and from the viewpoint of completeness of regular symplectic
reduction, and some developments around the work are give in
Wang and Zhang \cite{wazh12}, Ratiu and Wang \cite{rawa12},
Wang \cite{wa15a}.
In this section, as the application of the theoretical result,
we shall regard the rigid spacecraft-rotor system
with the control torque $u$ acting on the rotor
as a regular point reducible RCH system on the generalization
of rotation group $\textmd{SO}(3)\times \mathbb{R}$ and on the
generalization of Euclidean group $\textmd{SE}(3)\times \mathbb{R}$,
respectively. We shall give the $R_p$- reduced controlled spacecraft-rotor
systems in the two cases by calculation in detail,
We also follow the notations and conventions introduced
in  Marsden et al \cite{mawazh10}, Wang \cite{wa18} and Wang \cite{wa13d}.

\subsection{Spacecraft-Rotor System with Coincident Centers}

We first give the regular point reduction of the controlled rigid spacecraft-rotor system
with coincident centers of buoyancy and gravity. Assume that Lie
group $G=\textmd{SO}(3)$ acts freely and properly on
$Q=\textmd{SO}(3)\times \mathbb{R}$ by the left translation on the first
factor $\textmd{SO}(3)$, and the trivial action on the second factor
$\mathbb{R}$. By using the local left trivialization of $T^\ast \textmd{SO}(3)=
\textmd{SO}(3)\times \mathfrak{so}^\ast(3) $, then the action of
$\textmd{SO}(3)$ on phase space $T^\ast Q= T^\ast
\textmd{SO}(3)\times T^\ast \mathbb{R}\cong
\textmd{SO}(3)\times \mathfrak{so}^\ast(3)\times\mathbb{R}\times\mathbb{R}^{*}$
is by cotangent lift of left
translation on $\textmd{SO}(3)$ at the identity, that is, $\Phi^{T*}:
\textmd{SO}(3)\times T^\ast Q \cong
\textmd{SO}(3)\times \textmd{SO}(3)\times \mathfrak{so}^\ast(3)
\times \mathbb{R} \times \mathbb{R}^*
\to T^* Q \cong \textmd{SO}(3)\times
\mathfrak{so}^\ast(3)\times \mathbb{R} \times \mathbb{R}^*,$ given by
$\Phi^{T*}(B,(A,\Pi,\alpha,l))=(BA,\Pi,\alpha,l)$, for any $A,B\in
\textmd{SO}(3), \; \Pi \in \mathfrak{so}^\ast(3), \; \alpha \in \mathbb{R}, \; l \in
\mathbb{R}^*$. Assume that the action is free, proper and symplectic,
and the orbit space $(T^* Q)/ \textmd{SO}(3)$ is a smooth manifold
and $\pi: T^*Q \rightarrow (T^*Q )/ \textmd{SO}(3) $ is a smooth submersion. Since
$\textmd{SO}(3)$ acts trivially on $\mathfrak{so}^\ast(3)$ and
$\mathbb{R}\times \mathbb{R}^*$, it follows that $(T^\ast Q)/ \textmd{SO}(3)$ is
diffeomorphic to $\mathfrak{so}^\ast(3) \times \mathbb{R}
\times \mathbb{R}^*$.\\

We know that $\mathfrak{so}^\ast(3)$ is a Poisson manifold with
respect to its rigid body Lie-Poisson bracket defined by
\begin{equation}
\{F,K\}_{\mathfrak{so}^\ast(3)}(\Pi)= -\Pi\cdot (\nabla_\Pi
F\times \nabla_\Pi K), \;\; \forall F,K\in
C^\infty(\mathfrak{so}^\ast(3)),\;\; \Pi \in
\mathfrak{so}^\ast(3).\label{3.1}
\end{equation}
For $\mu \in \mathfrak{so}^\ast(3)$, the co-adjoint orbit
$\mathcal{O}_\mu \subset \mathfrak{so}^\ast(3)$ has the induced
orbit symplectic form $\omega^{-}_{\mathcal{O}_\mu}$, which
coincides with the restriction of the Lie-Poisson bracket on
$\mathfrak{so}^\ast(3)$ to the co-adjoint orbit $\mathcal{O}_\mu$.
From the Symplectic Stratification theorem we know that the
co-adjoint orbits $(\mathcal{O}_\mu, \omega_{\mathcal{O}_\mu}^{-}),
\; \mu\in \mathfrak{so}^\ast(3),$ form the symplectic leaves of the
Poisson manifold
$(\mathfrak{so}^\ast(3),\{\cdot,\cdot\}_{\mathfrak{so}^\ast(3)}). $
Let $\omega_{\mathbb{R}}$ be the canonical symplectic form on
$T^\ast \mathbb{R} \cong \mathbb{R} \times \mathbb{R}^*$, which is
given by
\begin{equation}
 \omega_{\mathbb{R}}((\theta_1, \lambda_1),(\theta_2,
\lambda_2))=<\lambda_2,\theta_1> -<\lambda_1,\theta_2>,
\label{3.2} \end{equation}
where $(\theta_i, \lambda_i)\in \mathbb{R}\times
\mathbb{R}^*, \; i=1,2$, $<\cdot,\cdot>$ is the standard inner product
on $\mathbb{R}\times \mathbb{R}^*$. It induces a canonical Poisson
bracket $\{\cdot,\cdot\}_{\mathbb{R}}$ on $T^\ast \mathbb{R}$, which
is given by
\begin{equation}
\{F,K\}_{\mathbb{R}}(\theta,\lambda)= \frac{\partial
F}{\partial \theta} \frac{\partial K}{\partial \lambda}-
\frac{\partial K}{\partial \theta}\frac{\partial F}{\partial
\lambda} .\label{3.3}
\end{equation}
See Marsden and Ratiu \cite{mara99}. Thus, we can induce a symplectic form
$\tilde{\omega}^{-}_{\mathcal{O}_\mu \times \mathbb{R} \times
\mathbb{R}^*}= \pi_{\mathcal{O}_\mu}^\ast
\omega^{-}_{\mathcal{O}_\mu}+ \pi_{\mathbb{R}}^\ast
\omega_{\mathbb{R}}$ on the smooth manifold $\mathcal{O}_\mu \times
\mathbb{R} \times \mathbb{R}^*$, where the maps
$\pi_{\mathcal{O}_\mu}: \mathcal{O}_\mu \times \mathbb{R} \times
\mathbb{R}^* \to \mathcal{O}_\mu$ and $\pi_{\mathbb{R}}:
\mathcal{O}_\mu \times \mathbb{R} \times \mathbb{R}^* \to \mathbb{R}
\times \mathbb{R}^*$ are canonical projections, and can induce a Poisson
bracket $\{\cdot,\cdot\}_{-}= \pi_{\mathfrak{so}^\ast(3)}^\ast
\{\cdot,\cdot\}_{\mathfrak{so}^\ast(3)}+ \pi_{\mathbb{R}}^\ast
\{\cdot,\cdot\}_{\mathbb{R}}$ on the smooth manifold
$\mathfrak{so}^\ast(3)\times \mathbb{R} \times \mathbb{R}^*$, where
the maps $\pi_{\mathfrak{so}^\ast(3)}: \mathfrak{so}^\ast(3) \times
\mathbb{R} \times \mathbb{R}^* \to \mathfrak{so}^\ast(3)$ and
$\pi_{\mathbb{R}}: \mathfrak{so}^\ast(3) \times \mathbb{R} \times
\mathbb{R}^* \to \mathbb{R} \times \mathbb{R}^*$ are canonical
projections, and such that $(\mathcal{O}_\mu \times \mathbb{R}\times
\mathbb{R}^*,\tilde{\omega}_{\mathcal{O}_\mu \times \mathbb{R} \times
\mathbb{R}^*}^{-})$ is a symplectic leaf of the Poisson manifold
$(\mathfrak{so}^\ast(3) \times \mathbb{R} \times \mathbb{R}^*, \{\cdot,\cdot\}_{-}). $\\

On the other hand, from $T^\ast Q = T^\ast \textmd{SO}(3) \times
T^\ast \mathbb{R}$ we know that there is a canonical symplectic form
$\omega_Q= \pi^\ast_{\textmd{SO}(3)} \omega_0 +\pi^\ast_{\mathbb{R}}
\omega_{\mathbb{R}}$ on $T^\ast Q$, where $\omega_0$ is the canonical
symplectic form on $T^\ast \textmd{SO}(3)$ and the maps
$\pi_{\textmd{SO}(3)}: Q= \textmd{SO}(3)\times \mathbb{R} \to
\textmd{SO}(3)$ and $\pi_{\mathbb{R}}: Q= \textmd{SO}(3)\times \mathbb{R} \to \mathbb{R}$
are canonical projections. Assume that the cotangent lift of left
$\textmd{SO}(3)$-action $\Phi^{T*}: \textmd{SO}(3) \times T^\ast Q \to
T^\ast Q$ is symplectic with respect to $\omega_Q$, and admits an associated
$\operatorname{Ad}^\ast$-equivariant momentum map $\mathbf{J}_Q:
T^\ast Q \to \mathfrak{so}^\ast(3)$ such that $\mathbf{J}_Q\cdot
\pi^\ast_{\textmd{SO}(3)}=\mathbf{J}_{\textmd{SO}(3)}$, where
$\mathbf{J}_{\textmd{SO}(3)}: T^\ast \textmd{SO}(3) \rightarrow
\mathfrak{so}^\ast(3)$ is a momentum map of left
$\textmd{SO}(3)$-action on $T^\ast \textmd{SO}(3)$ and we assume that it exists, and
$\pi^\ast_{\textmd{SO}(3)}: T^\ast \textmd{SO}(3) \to T^\ast Q$. If
$\mu\in\mathfrak{so}^\ast(3)$ is a regular value of $\mathbf{J}_Q$,
then $\mu\in\mathfrak{so}^\ast(3)$ is also a regular value of
$\mathbf{J}_{\textmd{SO}(3)}$ and $\mathbf{J}_Q^{-1}(\mu)\cong
\mathbf{J}_{\textmd{SO}(3)}^{-1}(\mu)\times \mathbb{R} \times
\mathbb{R}^*$. Denote by $\textmd{SO}(3)_\mu=\{g\in
\textmd{SO}(3)|\operatorname{Ad}_g^\ast \mu=\mu \}$ the isotropy
subgroup of co-adjoint $\textmd{SO}(3)$-action at the point
$\mu\in\mathfrak{so}^\ast(3)$. It follows that $\textmd{SO}(3)_\mu$
acts also freely and properly on $\mathbf{J}_Q^{-1}(\mu)$, the
$R_p$-reduced space $(T^\ast
Q)_\mu=\mathbf{J}_Q^{-1}(\mu)/\textmd{SO}(3)_\mu\cong (T^\ast
\textmd{SO}(3))_\mu \times \mathbb{R} \times \mathbb{R}^*$ of $(T^\ast
Q,\omega_Q)$ at $\mu$, is a symplectic manifold with symplectic form
$\omega_\mu$ uniquely characterized by the relation $\pi_\mu^\ast
\omega_\mu=i_\mu^\ast \omega_Q=i_\mu^\ast \pi^\ast_{\textmd{SO}(3)}
\omega_0 +i_\mu^\ast \pi^\ast_{\mathbb{R}} \omega_{\mathbb{R}}$, where the map
$i_\mu:\mathbf{J}_Q^{-1}(\mu)\rightarrow T^\ast Q$ is the inclusion
and $\pi_\mu:\mathbf{J}_Q^{-1}(\mu)\rightarrow (T^\ast Q)_\mu$ is
the projection. From Abraham and Marsden \cite{abma78}, we
know that $((T^\ast \textmd{SO}(3))_\mu,\omega_\mu)$ is
symplectically diffeomorphic to
$(\mathcal{O}_\mu,\omega_{\mathcal{O}_\mu}^{-})$, and hence we have
that $((T^\ast Q)_\mu,\omega_\mu)$ is symplectically diffeomorphic
to $(\mathcal{O}_\mu \times \mathbb{R}\times
\mathbb{R}^*,\tilde{\omega}_{\mathcal{O}_\mu \times \mathbb{R} \times
\mathbb{R}^*}^{-}). $\\

From the expression $(2.1)$ of the Hamiltonian, we know that
$H(A,\Pi,\alpha,l)$ is invariant under the cotangent lift of the left
$\textmd{SO}(3)$-action $\Phi^*: \textmd{SO}(3)\times T^\ast Q \to
T^\ast Q$. From the rigid body Lie-Poisson
bracket on $\mathfrak{so}^\ast(3)$ and the Poisson bracket on
$T^\ast \mathbb{R}$, we can get the Poisson bracket on
$\mathfrak{so}^\ast(3)\times\mathbb{R}\times\mathbb{R}^*$, that is,
for $F,K: \mathfrak{so}^\ast(3)\times\mathbb{R}\times\mathbb{R}^* \to
\mathbb{R}, $ we have that
\begin{align}
\{F,K\}_{-}(\Pi,\alpha,l)=-\Pi\cdot(\nabla_\Pi F\times \nabla_\Pi
K)+ \{F,K\}_{\mathbb{R}}(\alpha,l).
\label{3.4} \end{align}
See Krishnaprasad and Marsden \cite{krma87}.
Hence, the Hamiltonian vector field $X_H$ of rigid spacecraft-rotor system is given by
\begin{align*}
X_{H}(\Pi)& =\{\Pi,\; H\}_{-}=
-\Pi\cdot(\nabla_\Pi\Pi\times\nabla_\Pi
H)+ \{\Pi,\; H\}_{\mathbb{R}}\\
&= -\nabla_\Pi\Pi\cdot(\nabla_\Pi H\times \Pi)+ (\frac{\partial \Pi}{\partial \alpha}
\frac{\partial H}{\partial l}- \frac{\partial
H}{\partial \alpha}\frac{\partial \Pi}{\partial l})\\
& =(\Pi_1,\Pi_2,\Pi_3)\times (\frac{\Pi_1}{ \bar{I}_1},\;\;
\frac{\Pi_2}{ \bar{I}_2}, \;\; \frac{(\Pi_3-
l)}{\bar{I}_3}) \\
&= ( \frac{(\bar{I}_2-\bar{I}_3)\Pi_2\Pi_3-
\bar{I}_2\Pi_2l }{\bar{I}_2\bar{I}_3}, \;\;
\frac{(\bar{I}_3-\bar{I}_1)\Pi_3\Pi_1+
\bar{I}_1\Pi_1l}{\bar{I}_3\bar{I}_1}, \;\;
\frac{(\bar{I}_1-\bar{I}_2)\Pi_1\Pi_2}{\bar{I}_1\bar{I}_2} ),
\end{align*}
since $\nabla_{\Pi_i}\Pi_i=1,\; \nabla_{\Pi_i}\Pi_j=0, \; i\neq j , \; i, j=1,2,3 $,
and $\nabla_{\Pi_k} H= \Pi_k / \bar{I}_k, \; k=1,2 $, $\nabla_{\Pi_3} H= (\Pi_3-l)/ \bar{I}_3, \;  \frac{\partial
\Pi}{\partial \alpha}= \frac{\partial H}{\partial
\alpha}=0. $

\begin{align*}
X_{H}(\alpha)& =\{\alpha,\; H\}_{-}=
-\Pi\cdot(\nabla_\Pi\alpha\times\nabla_\Pi
H)+ \{\alpha,\; H\}_{\mathbb{R}}\\
&= -\nabla_\Pi \alpha\cdot(\nabla_\Pi H\times \Pi)+ (\frac{\partial \alpha}{\partial \alpha}
\frac{\partial H}{\partial l}- \frac{\partial
H}{\partial \alpha}\frac{\partial \alpha}{\partial l})
= -\frac{(\Pi_3- l)}{\bar{I}_3} +\frac{l}{J_3},
\end{align*}
since $\nabla_{\Pi_i}\alpha=0, \;  i= 1,2,3$, $\frac{\partial \alpha}{\partial
\alpha}= 1, \; \frac{\partial H}{\partial \alpha}=0, $
and $\frac{\partial H}{\partial l}= -(\Pi_3-l)/\bar{I}_3
+\frac{l}{J_3}. $

\begin{align*}
X_{H}(l)& =\{l,\; H\}_{-}=
-\Pi\cdot(\nabla_\Pi l\times\nabla_\Pi
H)+ \{l,\; H\}_{\mathbb{R}}\\
&= -\nabla_\Pi l\cdot(\nabla_\Pi H\times \Pi)+ (\frac{\partial l}{\partial \alpha}
\frac{\partial H}{\partial l}- \frac{\partial
H}{\partial \alpha}\frac{\partial l}{\partial l})
=0,
\end{align*}
since $\nabla_{\Pi_i} l=0, \; i=1,2,3$ and $\frac{\partial l}{\partial \alpha}=
\frac{\partial H}{\partial \alpha}=0. $\\

Moreover, if we consider the rigid spacecraft-rotor system with a
control torque $u: T^\ast Q \to W $ acting on the rotor,
where the control subset $W\subset T^* Q $ is a fiber submanifold,
and assume that $u\in W $ is invariant under the cotangent lift $\Phi^{T^*}$ of left
$\textmd{SO}(3)$-action, and
the dynamical vector field of the regular point reducible
controlled rigid spacecraft-rotor system $(T^\ast Q,\textmd{SO}(3),\omega_Q,H,u)$
can be expressed by
\begin{align}
\tilde{X}= X_{(T^\ast Q,\textmd{SO}(3),\omega_Q,H,u)}= X_H+ \textnormal{vlift}(u),
\label{3.5} \end{align}
where $\textnormal{vlift}(u)= \textnormal{vlift}(u)\cdot X_H $ is
the change of $X_H$ under the action of the control torque $u$.
From the above expression of the dynamical vector
field of the controlled spacecraft-rotor system $(T^\ast Q,\textmd{SO}(3),\omega_Q,H,u)$,
we know that under the actions of the control torque $u$, in general, the dynamical vector
field is not Hamiltonian, and hence the regular point reducible
controlled rigid spacecraft-rotor system is not
yet a Hamiltonian system. However,
it is a dynamical system closed relative to a
Hamiltonian system, and it can be explored and studied by extending
the methods for the control torque $u$
in the study of the Marsden-Weinstein reducible Hamiltonian system
$(T^\ast Q,\textmd{SO}(3),\omega_Q,H)$,
see Marsden et al \cite{mawazh10} and Wang \cite{wa18}. \\

Since the Hamiltonian $H(A,\Pi,\alpha,l)$ is invariant under the cotangent lift $\Phi^{T^*}$ of the left
$\textmd{SO}(3)$-action, for the point $\Pi_0 =\mu \in \mathfrak{so}^\ast(3)$ is
the regular value of $\mathbf{J}_Q$, we have the $R_p$-reduced Hamiltonian
$h_\mu(\Pi, \alpha,l):\mathcal{O}_\mu
\times\mathbb{R}\times\mathbb{R}^*(\subset \mathfrak{so}^\ast
(3)\times\mathbb{R}\times\mathbb{R}^*)\to \mathbb{R}$ given by
$h_\mu(\Pi,\alpha,l)\cdot \pi_\mu
=H(A,\Pi,\alpha,l)|_{\mathcal{O}_\mu
\times\mathbb{R}\times\mathbb{R}^*}.$
Moreover, for the $R_p$-reduced
Hamiltonian $h_\mu(\Pi,\alpha,l): \mathcal{O}_\mu
\times\mathbb{R}\times\mathbb{R}^{*} \to \mathbb{R}$, we have the
Hamiltonian vector field
$X_{h_\mu}(K_\mu)=\{K_\mu,h_\mu\}_{-}|_{\mathcal{O}_\mu
\times\mathbb{R}\times\mathbb{R}^{*} }, $ where
$K_\mu(\Pi,\alpha,l): \mathcal{O}_\mu
\times\mathbb{R}\times\mathbb{R}^{*} \to \mathbb{R}.$
Assume that $u\in
W \cap \mathbf{J}^{-1}_Q(\mu)$ and the $R_p$-reduced control torque $u_\mu:
\mathcal{O}_\mu \times\mathbb{R}\times\mathbb{R}^{*} \to W_\mu
(\subset \mathcal{O}_\mu \times\mathbb{R}\times\mathbb{R}^{*}) $ is
given by $u_\mu(\Pi,\alpha,l)\cdot \pi_\mu=
u(A,\Pi,\alpha,l)|_{\mathcal{O}_\mu
\times\mathbb{R}\times\mathbb{R}^{*} }, $ where $\pi_\mu:
\mathbf{J}_Q^{-1}(\mu) \rightarrow \mathcal{O}_\mu
\times\mathbb{R}\times\mathbb{R}^{*}, \; W_\mu= \pi_\mu(W\cap \mathbf{J}^{-1}_Q(\mu)). $
The $R_p$-reduced controlled rigid spacecraft-rotor
system is the 4-tuple $(\mathcal{O}_\mu \times \mathbb{R} \times
\mathbb{R}^{*},\tilde{\omega}_{\mathcal{O}_\mu \times \mathbb{R}
\times \mathbb{R}^{*}}^{-},h_\mu,u_\mu), $ where
$\tilde{\omega}_{\mathcal{O}_\mu \times \mathbb{R} \times
\mathbb{R}^{*}}^{-}$ is the induced symplectic form from the Poisson bracket
on $\mathfrak{so}^\ast(3) \times \mathbb{R}\times \mathbb{R}^{*} ,$ such that
Hamiltonian vector field
$X_{h_\mu}(K_\mu)=\tilde{\omega}_{\mathcal{O}_\mu \times \mathbb{R}\times
\mathbb{R}^{*}}^{-}(X_{K_\mu}, X_{h_\mu})
=\{K_\mu,h_\mu\}_{-}|_{\mathcal{O}_\mu
\times\mathbb{R}\times\mathbb{R}^{*} },$
since $(\mathcal{O}_\mu \times \mathbb{R}\times
\mathbb{R}^*,\tilde{\omega}_{\mathcal{O}_\mu \times \mathbb{R} \times
\mathbb{R}^*}^{-})$ is a symplectic leaf of the Poisson manifold
$(\mathfrak{so}^\ast(3) \times \mathbb{R} \times \mathbb{R}^*, \{\cdot,\cdot\}_{-}). $
Moreover, assume that the dynamical vector field of the $R_p$-reduced controlled
spacecraft-rotor system $(\mathcal{O}_\mu \times \mathbb{R}
\times \mathbb{R}^*,\tilde{\omega}_{\mathcal{O}_\mu \times \mathbb{R}
\times \mathbb{R}^*}^{-},h_\mu,u_\mu)$ is expressed by
\begin{align} X_{(\mathcal{O}_\mu \times \mathbb{R} \times
\mathbb{R}^*,\tilde{\omega}_{\mathcal{O}_\mu \times \mathbb{R} \times
\mathbb{R}^*}^{-},h_\mu,u_\mu)} = X_{h_\mu} + \mbox{vlift}(u_\mu) ,
\label {3.6}\end{align}
where $\mbox{vlift}(u_\mu)= \mbox{vlift}(u_\mu)X_{h_\mu} \in
T(\mathcal{O}_\mu \times\mathbb{R}\times\mathbb{R}^*), $
is the change of $X_{h_\mu}$ under the action of the $R_p$-reduced control torque $u_\mu$.
The dynamical vector fields of the controlled
spacecraft-rotor system and the $R_p$-reduced controlled
spacecraft-rotor system satisfy the condition
\begin{equation}X_{(\mathcal{O}_\mu \times \mathbb{R} \times
\mathbb{R}^*,\tilde{\omega}_{\mathcal{O}_\mu \times \mathbb{R} \times
\mathbb{R}^*}^{-},h_\mu,u_\mu)}\cdot \pi_\mu=T\pi_\mu\cdot X_{(T^\ast
Q,\textmd{SO}(3),\omega_Q,H,u)}\cdot i_\mu. \label{3.7}\end{equation}
See Marsden et al \cite{mawazh10} and Wang \cite{wa18}.\\

To sum up the above discussion, we have the following theorem.
\begin{theo}
In the case of coincident centers of buoyancy and gravity, the rigid
spacecraft-rotor system with the control torque $u$ acting on the
rotor, that is, the 5-tuple $(T^\ast Q, \textmd{SO}(3), \omega_Q, H,
u ), $ where $Q= \textmd{SO}(3)\times \mathbb{R}, $ is a regular point
reducible RCH system. For a point $\mu \in \mathfrak{so}^\ast(3)$,
the regular value of the momentum map $\mathbf{J}_Q:
\textmd{SO}(3)\times \mathfrak{so}^\ast(3) \times \mathbb{R} \times
\mathbb{R}^* \to \mathfrak{so}^\ast(3)$, the $R_p$-reduced
controlled rigid spacecraft-rotor system is the 4-tuple $(\mathcal{O}_\mu \times \mathbb{R} \times
\mathbb{R}^*,\tilde{\omega}_{\mathcal{O}_\mu \times \mathbb{R} \times
\mathbb{R}^*}^{-},h_\mu,u_\mu), $ where $\mathcal{O}_\mu \subset
\mathfrak{so}^\ast(3)$ is the co-adjoint orbit,
$\tilde{\omega}_{\mathcal{O}_\mu \times \mathbb{R} \times
\mathbb{R}^*}^{-}$ is the induced symplectic form
on $\mathcal{O}_\mu \times
\mathbb{R}\times \mathbb{R}^* $,
$h_\mu(\Pi,\alpha,l)\cdot \pi_\mu=H(A,\Pi,\alpha,l)|_{\mathcal{O}_\mu
\times\mathbb{R}\times\mathbb{R}^*}$, $u_\mu(\Pi,\alpha,l)\cdot
\pi_\mu= u(A,\Pi,\alpha,l)|_{\mathcal{O}_\mu
\times\mathbb{R}\times\mathbb{R}^*}$, and
the dynamical vector field of the $R_p$-reduced controlled
spacecraft-rotor system $(\mathcal{O}_\mu \times \mathbb{R}
\times \mathbb{R}^*,\tilde{\omega}_{\mathcal{O}_\mu \times \mathbb{R}
\times \mathbb{R}^*}^{-},h_\mu,u_\mu)$ satisfies (3.6) and (3.7).
\end{theo}

\subsection{Spacecraft-Rotor System with Non-coincident Centers}

In the following we shall give the regular point reduction of the controlled rigid
spacecraft-rotor system with non-coincident centers of buoyancy and
gravity. Because the drift in the direction of gravity breaks the
symmetry and the spacecraft-rotor system is no longer
$\textmd{SO}(3)$ invariant. In this case, its physical phase space
is $T^\ast \textmd{SO}(3)\times T^* S^1$ and the symmetry group is
$S^1$, regarded as rotations about the third principal axis, that
is, the axis of gravity. By the semidirect product reduction
theorem, see Marsden et al \cite{mamiorpera07}, we know that the
reduction of $T^\ast \textmd{SO}(3)$ by $S^1$ gives a space which is
symplectically diffeomorphic to the reduced space obtained by the
reduction of $T^\ast \textmd{SE}(3)$ by left action of
$\textmd{SE}(3)$, that is the coadjoint orbit $\mathcal{O}_{(\mu,a)}
\subset \mathfrak{se}^\ast(3)\cong T^\ast
\textmd{SE}(3)/\textmd{SE}(3)$. In fact, in this case, we can
identify the phase space $T^\ast \textmd{SO}(3)$ with the reduction
of the cotangent bundle of the special Euclidean group
$\textmd{SE}(3)=\textmd{SO}(3)\circledS \mathbb{R}^3$ by the
Euclidean translation subgroup $\mathbb{R}^3$ and identifies the
symmetry group $S^1$ with isotropy group $G_a=\{ A\in
\textmd{SO}(3)\mid Aa=a \}=S^1$, which is Abelian and
$(G_a)_{\mu_a}= G_a =S^1,\; \forall \mu_a \in \mathfrak{g}^\ast_a$,
where $a$ is a vector aligned with the direction of gravity and
where $\textmd{SO}(3)$ acts on $\mathbb{R}^3$ in the standard way.\\

Assume that Lie group $G=\textmd{SE}(3)$ acts freely and properly on
$Q=\textmd{SE}(3)\times \mathbb{R}$ by the left translation on the first
factor $\textmd{SE}(3)$, and the trivial action on the second factor
$\mathbb{R}$. By using the local left trivialization of $T^\ast \textmd{SE}(3)=
\textmd{SE}(3)\times \mathfrak{se}^\ast(3) $, then the action of
$\textmd{SE}(3)$ on phase space $T^\ast Q= T^\ast
\textmd{SE}(3)\times T^\ast \mathbb{R}$ is by cotangent lift of left
translation on $\textmd{SE}(3)$ at the identity, that is, $\Phi^{T*}:
\textmd{SE}(3)\times T^\ast Q \cong
\textmd{SE}(3)\times \textmd{SE}(3)\times \mathfrak{se}^\ast(3)
\times \mathbb{R} \times \mathbb{R}^*\to T^* Q \cong \textmd{SE}(3)\times
\mathfrak{se}^\ast(3)\times \mathbb{R} \times \mathbb{R}^*,$ given by
$\Phi^{T*}((B,b)(A,c,\Pi,\Gamma,\alpha,l))=(BA,b+Bc,\Pi,\Gamma,\alpha,l)$,
for any $A,B\in \textmd{SO}(3), \; \Pi \in \mathfrak{so}^\ast(3), \;
b, c, \Gamma \in \mathbb{R}^3,  \; (\Pi,\Gamma) \in
\mathfrak{se}^\ast(3) \; \alpha \in \mathbb{R}, \; l \in \mathbb{R}^*$. Assume that
the action is free, proper and symplectic, and the orbit space $(T^* Q)/ \textmd{SE}(3)$
is a smooth manifold and $\pi: T^*Q \rightarrow (T^*Q )/
\textmd{SE}(3) $ is a smooth submersion. Since $\textmd{SE}(3)$ acts
trivially on $\mathfrak{se}^\ast(3)$ and $\mathbb{R}\times \mathbb{R}^*$, it follows
that $(T^\ast Q)/ \textmd{SE}(3)$ is diffeomorphic to
$\mathfrak{se}^\ast(3) \times \mathbb{R}
\times \mathbb{R}^*$.\\

We know that $\mathfrak{se}^\ast(3)$ is a Poisson manifold with
respect to its heavy top Lie-Poisson bracket defined by
\begin{equation}   \{F,K\}_{\mathfrak{se}^\ast(3)}(\Pi,\Gamma)= -\Pi\cdot (\nabla_\Pi
F\times \nabla_\Pi K)-\Gamma\cdot(\nabla_\Pi F\times \nabla_\Gamma
K-\nabla_\Pi K\times \nabla_\Gamma F), \;\; \label{3.8}
\end{equation}
where $F, K \in C^\infty(\mathfrak{se}^\ast(3)), \; (\Pi,\Gamma) \in
\mathfrak{se}^\ast(3)$, see Marsden et al.
\cite{mamiorpera07}. For $(\mu,a) \in \mathfrak{se}^\ast(3)$, the
co-adjoint orbit $\mathcal{O}_{(\mu,a)} \subset
\mathfrak{se}^\ast(3)$ has the induced orbit symplectic form
$\omega^{-}_{\mathcal{O}_{(\mu,a)}}$, which coincides with the
restriction of the Lie-Poisson bracket on $\mathfrak{se}^\ast(3)$ to
the co-adjoint orbit $\mathcal{O}_{(\mu,a)}$.
From the Symplectic Stratification theorem we know that the co-adjoint
orbits $(\mathcal{O}_{(\mu,a)}, \omega_{\mathcal{O}_{(\mu,a)}}^{-})$,
$ (\mu,a)\in \mathfrak{se}^\ast(3),$ form the symplectic leaves of
the Poisson manifold
$(\mathfrak{se}^\ast(3),\{\cdot,\cdot\}_{\mathfrak{se}^\ast(3)}). $
Let $\omega_{\mathbb{R}}$ be the canonical symplectic form on
$T^\ast \mathbb{R} \cong \mathbb{R} \times \mathbb{R}$ given by
$(3.2)$, and it induces a canonical Poisson bracket
$\{\cdot,\cdot\}_{\mathbb{R}}$ on $T^\ast \mathbb{R}$ given by
$(3.3)$. Thus, we can induce a symplectic form
$\tilde{\omega}^{-}_{\mathcal{O}_{(\mu,a)} \times \mathbb{R} \times
\mathbb{R}}= \pi_{\mathcal{O}_{(\mu,a)}}^\ast
\omega^{-}_{\mathcal{O}_{(\mu,a)}}+ \pi_{\mathbb{R}}^\ast
\omega_{\mathbb{R}}$ on the smooth manifold $\mathcal{O}_{(\mu,a)}
\times \mathbb{R} \times \mathbb{R}$, where the maps
$\pi_{\mathcal{O}_{(\mu,a)}}: \mathcal{O}_{(\mu,a)} \times
\mathbb{R} \times \mathbb{R} \to \mathcal{O}_{(\mu,a)}$ and
$\pi_{\mathbb{R}}: \mathcal{O}_{(\mu,a)} \times \mathbb{R} \times
\mathbb{R} \to \mathbb{R} \times \mathbb{R}$ are canonical
projections, and can induce a Poisson bracket $\{\cdot,\cdot\}_{-}=
\pi_{\mathfrak{se}^\ast(3)}^\ast
\{\cdot,\cdot\}_{\mathfrak{se}^\ast(3)}+ \pi_{\mathbb{R}}^\ast
\{\cdot,\cdot\}_{\mathbb{R}}$ on the smooth manifold
$\mathfrak{se}^\ast(3)\times \mathbb{R} \times \mathbb{R}$, where
the maps $\pi_{\mathfrak{se}^\ast(3)}: \mathfrak{se}^\ast(3) \times
\mathbb{R} \times \mathbb{R} \to \mathfrak{se}^\ast(3)$ and
$\pi_{\mathbb{R}}: \mathfrak{se}^\ast(3) \times \mathbb{R} \times
\mathbb{R} \to \mathbb{R} \times \mathbb{R}$ are canonical
projections, and such that $(\mathcal{O}_{(\mu,a)} \times
\mathbb{R}\times \mathbb{R},\tilde{\omega}_{\mathcal{O}_{(\mu,a)}
\times \mathbb{R} \times \mathbb{R}}^{-})$ is a symplectic leaf of
the Poisson manifold
$(\mathfrak{se}^\ast(3) \times \mathbb{R} \times \mathbb{R}, \{\cdot,\cdot\}_{-}). $\\

On the other hand, from $T^\ast Q = T^\ast \textmd{SE}(3) \times
T^\ast \mathbb{R}$ we know that there is a canonical symplectic form
$\omega_Q= \pi^\ast_{\textmd{SE}(3)} \omega_1 +\pi^\ast_{\mathbb{R}}
\omega_{\mathbb{R}}$ on $T^\ast Q$, where $\omega_1$ is the canonical
symplectic form on $T^\ast \textmd{SE}(3)$ and the maps
$\pi_{\textmd{SE}(3)}: Q= \textmd{SE}(3)\times \mathbb{R} \to
\textmd{SE}(3)$ and $\pi_{\mathbb{R}}: Q= \textmd{SE}(3)\times \mathbb{R} \to \mathbb{R}$
are canonical projections. Assume that the cotangent lift of left
$\textmd{SE}(3)$-action $\Phi^{T*}: \textmd{SE}(3) \times T^\ast Q \to
T^\ast Q$ is symplectic with respect to $\omega_Q$, and admits an associated
$\operatorname{Ad}^\ast$-equivariant momentum map $\mathbf{J}_Q:
T^\ast Q \to \mathfrak{se}^\ast(3)$ such that $\mathbf{J}_Q\cdot
\pi^\ast_{\textmd{SE}(3)}=\mathbf{J}_{\textmd{SE}(3)}$, where
$\mathbf{J}_{\textmd{SE}(3)}: T^\ast \textmd{SE}(3) \rightarrow
\mathfrak{se}^\ast(3)$ is a momentum map of left
$\textmd{SE}(3)$-action on $T^\ast \textmd{SE}(3)$ and assume that it exists, and
$\pi^\ast_{\textmd{SE}(3)}: T^\ast \textmd{SE}(3) \to T^\ast Q$. If
$(\mu,a)\in\mathfrak{se}^\ast(3)$ is a regular value of
$\mathbf{J}_Q$, then $(\mu,a)\in\mathfrak{se}^\ast(3)$ is also a
regular value of $\mathbf{J}_{\textmd{SE}(3)}$ and
$\mathbf{J}_Q^{-1}(\mu,a)\cong
\mathbf{J}_{\textmd{SE}(3)}^{-1}(\mu,a)\times \mathbb{R} \times
\mathbb{R}^*$. Denote by $\textmd{SE}(3)_{(\mu,a)}=\{g\in
\textmd{SE}(3)|\operatorname{Ad}_g^\ast (\mu,a)=(\mu,a) \}$ the
isotropy subgroup of co-adjoint $\textmd{SE}(3)$-action at the point
$(\mu,a)\in\mathfrak{se}^\ast(3)$. It follows that
$\textmd{SE}(3)_{(\mu,a)}$ acts also freely and properly on
$\mathbf{J}_Q^{-1}(\mu,a)$, the $R_p$-reduced space $(T^\ast
Q)_{(\mu,a)}=\mathbf{J}_Q^{-1}(\mu,a)/\textmd{SE}(3)_{(\mu,a)}\cong
(T^\ast \textmd{SE}(3))_{(\mu,a)} \times \mathbb{R} \times
\mathbb{R}^*$ of $(T^\ast Q,\omega_Q)$ at $(\mu,a)$, is a symplectic
manifold with symplectic form $\omega_{(\mu,a)}$ uniquely
characterized by the relation $\pi_{(\mu,a)}^\ast
\omega_{(\mu,a)}=i_{(\mu,a)}^\ast \omega_Q=i_{(\mu,a)}^\ast
\pi^\ast_{\textmd{SE}(3)} \omega_1 +i_{(\mu,a)}^\ast \pi^\ast_{\mathbb{R}}
\omega_{\mathbb{R}}$, where the map
$i_{(\mu,a)}:\mathbf{J}_Q^{-1}(\mu,a)\rightarrow T^\ast Q$ is the
inclusion and $\pi_{(\mu,a)}:\mathbf{J}_Q^{-1}(\mu,a)\rightarrow
(T^\ast Q)_{(\mu,a)}$ is the projection. From Abraham and
Marsden \cite{abma78}, we know that $((T^\ast
\textmd{SE}(3))_{(\mu,a)},\omega_{(\mu,a)})$ is symplectically
diffeomorphic to
$(\mathcal{O}_{(\mu,a)},\omega_{\mathcal{O}_{(\mu,a)}}^{-})$, and
hence we have that $((T^\ast Q)_{(\mu,a)},\omega_{(\mu,a)})$ is
symplectically diffeomorphic to $(\mathcal{O}_{(\mu,a)} \times
\mathbb{R}\times \mathbb{R}^*,\tilde{\omega}_{\mathcal{O}_{(\mu,a)}
\times \mathbb{R} \times \mathbb{R}^*}^{-}). $\\

From the expression $(2.2)$ of the Hamiltonian, we know that
$H(A,c,\Pi,\Gamma,\alpha,l)$ is invariant under the cotangent lift of the left
$\textmd{SE}(3)$-action $\Phi^{T*}: \textmd{SE}(3)\times T^\ast Q \to
T^\ast Q$. Moreover, from the heavy top Lie-Poisson
bracket on $\mathfrak{se}^\ast(3)$ and the Poisson bracket on
$T^\ast \mathbb{R}$, we can get the Poisson bracket
on $\mathfrak{se}^\ast(3)\times \mathbb{R}\times \mathbb{R}^{*} $,
that is, for $F,K: \mathfrak{se}^\ast(3)\times \mathbb{R}\times
\mathbb{R}^{*} \to \mathbb{R}, $ we have that
\begin{align} \{F,K\}_{-}(\Pi,\Gamma,\alpha,l)
= -\Pi\cdot(\nabla_\Pi F\times\nabla_\Pi K)-\Gamma\cdot(\nabla_\Pi
F\times \nabla_\Gamma K-\nabla_\Pi K\times \nabla_\Gamma F)+
\{F,K\}_{\mathbb{R}}(\alpha,l).
\label{3.9} \end{align}
Hence, the Hamiltonian vector field of the rigid spacecraft-rotor system is given by
\begin{align*}
& X_{H}(\Pi) = \{\Pi,\; H \}_{-}= -\Pi\cdot(\nabla_\Pi\Pi\times\nabla_\Pi
H) -\Gamma\cdot(\nabla_\Pi\Pi\times\nabla_\Gamma
H-\nabla_\Pi H \times\nabla_\Gamma\Pi) + \{\Pi,\; H \}_{\mathbb{R}}\\
& =(\Pi_1,\Pi_2,\Pi_3)\times (\frac{\Pi_1}{ \bar{I}_1}, \frac{\Pi_2}{
\bar{I}_2}, \frac{(\Pi_3-l)}{\bar{I}_3})
+gh(\Gamma_1,\Gamma_2,\Gamma_3)\times (\chi_1,\chi_2,\chi_3)
+ (\frac{\partial \Pi}{\partial \alpha}
\frac{\partial H}{\partial l}- \frac{\partial
H}{\partial \alpha}\frac{\partial \Pi}{\partial l})\\
& = ( \frac{(\bar{I}_2-\bar{I}_3)\Pi_2\Pi_3-
\bar{I}_2\Pi_2l}{\bar{I}_2\bar{I}_3}+
gh(\Gamma_2\chi_3-\Gamma_3\chi_2), \;\;
\frac{(\bar{I}_3-\bar{I}_1)\Pi_3\Pi_1+
\bar{I}_1\Pi_1l}{\bar{I}_3\bar{I}_1}+gh(\Gamma_3\chi_1-\Gamma_1\chi_3), \\
& \;\;\;\;\;\;
\frac{(\bar{I}_1-\bar{I}_2)\Pi_1\Pi_2}{\bar{I}_1\bar{I}_2}
 + gh(\Gamma_1\chi_2-\Gamma_2\chi_1) ),
\end{align*}
since $\nabla_{\Pi_i}\Pi_i=1,\; \nabla_{\Pi_i}\Pi_j=0, \; i\neq j, \;  \nabla_{\Pi_i}\Gamma_j=\nabla_{\Gamma_i}\Pi_j=0, $
and $\chi=(\chi_1,\chi_2,\chi_3), \; \nabla_{\Gamma_j}
H= gh\chi_j, \; \nabla_{\Pi_3} H= (\Pi_3-l )/\bar{I}_3 ,\; \nabla_{\Pi_k}
H= \Pi_k/\bar{I}_k, \; \frac{\partial
\Pi}{\partial \alpha}= \frac{\partial H}{\partial
\alpha}=0, \; i, j=1,2,3, \; k=1,2. $

\begin{align*}
 X_{H}(\Gamma)&
= \{\Gamma,\; H\}_{-} =-\Pi\cdot(\nabla_\Pi\Gamma\times\nabla_\Pi
H)-\Gamma\cdot(\nabla_\Pi\Gamma\times\nabla_\Gamma
H-\nabla_\Pi H \times\nabla_\Gamma\Gamma)+  \{\Gamma,\; H\}_{\mathbb{R}} \\
& =\nabla_\Gamma\Gamma\cdot (\Gamma\times\nabla_\Pi
H)+ (\frac{\partial \Gamma}{\partial \alpha}
\frac{\partial H}{\partial l}- \frac{\partial
H}{\partial \alpha}\frac{\partial \Gamma}{\partial l})
= (\Gamma_1,\Gamma_2,\Gamma_3)\times (\frac{\Pi_1}{ \bar{I}_1},\;\;
\frac{\Pi_2}{ \bar{I}_2}, \;\; \frac{(\Pi_3-l)}{\bar{I}_3}) \\
&= ( \frac{\bar{I}_2\Gamma_2\Pi_3-\bar{I}_3\Gamma_3\Pi_2-
\bar{I}_2\Gamma_2l}{\bar{I}_2\bar{I}_3}, \;\;
\frac{\bar{I}_3\Gamma_3\Pi_1-\bar{I}_1\Gamma_1\Pi_3+
\bar{I}_1\Gamma_1l}{\bar{I}_3\bar{I}_1},\;\;
\frac{\bar{I}_1\Gamma_1\Pi_2-\bar{I}_2\Gamma_2\Pi_1}{\bar{I}_1\bar{I}_2} ),
\end{align*}
since $\nabla_{\Gamma_i}\Gamma_i =1, \; \nabla_{\Gamma_i}\Gamma_j=0, \; i\neq j,
\; \nabla_{\Pi_i}\Gamma_j =0,  \; \frac{\partial
\Gamma_j}{\partial \alpha}= \frac{\partial H}{\partial
\alpha}=0, \; i, j= 1,2,3, $ and
$\nabla_{\Pi_3} H= (\Pi_3-l )/\bar{I}_3 ,\; \nabla_{\Pi_k}
H= \Pi_k/\bar{I}_k , \;  k=1,2.$

\begin{align*}
 X_{H}(\alpha)& = \{\alpha,\; H \}_{-}= -\Pi\cdot(\nabla_\Pi \alpha \times\nabla_\Pi
H) -\Gamma\cdot(\nabla_\Pi \alpha\times\nabla_\Gamma
H-\nabla_\Pi H \times\nabla_\Gamma \alpha) + \{\alpha,\; H \}_{\mathbb{R}}\\
& =(\frac{\partial \alpha}{\partial \alpha}
\frac{\partial H}{\partial l}- \frac{\partial
H}{\partial \alpha}\frac{\partial \alpha}{\partial l})=  -\frac{(\Pi_3- l)}{ \bar{I}_3}
+\frac{l}{J_3},
\end{align*}
since $\nabla_{\Pi_i}\alpha=\nabla_{\Gamma_i}\alpha=0, \; i= 1,2,3$,
$\frac{\partial \alpha}{\partial \alpha}= 1, \; \frac{\partial H}{\partial \alpha}=0, $
and $\frac{\partial H}{\partial l}= -(\Pi_3-l)/\bar{I}_3
+\frac{l}{J_3}. $

\begin{align*}
 X_{H}(l)& = \{l,\; H \}_{-}= -\Pi\cdot(\nabla_\Pi l \times\nabla_\Pi
H) -\Gamma\cdot(\nabla_\Pi l \times\nabla_\Gamma
H-\nabla_\Pi H \times\nabla_\Gamma l) + \{l,\; H \}_{\mathbb{R}}\\
& = (\frac{\partial l}{\partial \alpha}
\frac{\partial H}{\partial l}- \frac{\partial
H}{\partial \alpha}\frac{\partial l}{\partial l})=0,
\end{align*}
since $\nabla_{\Pi_i} l=\nabla_{\Gamma_i} l=0, \; i=1,2,3,$ and $\frac{\partial l}{\partial \alpha}=
\frac{\partial H}{\partial \alpha}=0. $\\

Moreover, if we consider the rigid spacecraft-rotor system with a control torque $u: T^\ast
Q \to W $ acting on the rotors, where the control subset $W\subset T^* Q $ is a fiber submanifold,
and assume that $u\in W $ is invariant under the cotangent lift
$\Phi^{T*}$ of the left $\textmd{SE}(3)$-action, and
the dynamical vector field of the regular point reducible
controlled spacecraft-rotor system $(T^\ast Q,\textmd{SE}(3),\omega_Q,H,u)$
can be expressed by
\begin{align}
\tilde{X}= X_{(T^\ast Q,\textmd{SE}(3),\omega_Q,H,u)}= X_H+ \textnormal{vlift}(u),
\label{3.10} \end{align}
where $\textnormal{vlift}(u)= \textnormal{vlift}(u)\cdot X_H $
is the change of $X_H$ under the action of the control torque $u$.
From the above expression of the dynamical vector
field of the spacecraft-rotor system $(T^\ast Q,\textmd{SE}(3),\omega_Q,H,u)$,
we know that under the actions of the control torque $u$, in general, the dynamical vector
field is not Hamiltonian, and hence the regular point reducible
controlled rigid spacecraft-rotor system is not
yet a Hamiltonian system. However,
it is a dynamical system closed relative to a
Hamiltonian system, and it can be explored and studied by extending
the methods for the control torque $u$
in the study of the Marsden-Weinstein reducible Hamiltonian system
$(T^\ast Q,\textmd{SE}(3),\omega_Q,H)$,
see Marsden et al \cite{mawazh10} and Wang \cite{wa18}. \\

Since the Hamiltonian
$H(A,c,\Pi,\Gamma,\alpha,l)$ is invariant under the cotangent lift $\Phi^{T*}$ of the left
$\textmd{SE}(3)$-action, for the point $(\Pi_0,\Gamma_0)=(\mu,a)\in
\mathfrak{se}^\ast(3)$ is the regular value of $\mathbf{J}_Q$, we
have the $R_p$-reduced Hamiltonian
$h_{(\mu,a)}(\Pi,\Gamma, \alpha,l): \mathcal{O}_{(\mu,a)}
\times\mathbb{R}\times\mathbb{R}^*(\subset \mathfrak{se}^\ast
(3)\times\mathbb{R}\times\mathbb{R}^*)\to \mathbb{R}$ given by
$h_{(\mu,a)}(\Pi,\Gamma,\alpha,l)\cdot
\pi_{(\mu,a)}= H(A,c,\Pi,\Gamma,\alpha,l)|_{\mathcal{O}_{(\mu,a)}
\times\mathbb{R}\times\mathbb{R}^*}.$
Moreover, for the $R_p$-reduced
Hamiltonian $h_{(\mu,a)}(\Pi,\Gamma, \alpha, l): \mathcal{O}_{(\mu,a)}\times
\mathbb{R}\times \mathbb{R}^{*} \to \mathbb{R}$, we have the
Hamiltonian vector field
$$X_{h_{(\mu,a)}}(K_{(\mu,a)})=\{K_{(\mu,a)},h_{(\mu,a)}\}_{-}|_{\mathcal{O}_{(\mu,a)}
\times \mathbb{R}\times \mathbb{R}^{*}}, $$ where
$K_{(\mu,a)}(\Pi,\Gamma, \alpha, l): \mathcal{O}_{(\mu,a)}\times
\mathbb{R}\times \mathbb{R}^{*} \to \mathbb{R}.$
Assume that $u\in W \cap
\mathbf{J}^{-1}_Q(\mu, a )$ and the  $R_p$-reduced control torque $u_{(\mu,a)}:
\mathcal{O}_{(\mu,a)} \times\mathbb{R}\times\mathbb{R}^{*} \to W_{(\mu,a)}
(\subset \mathcal{O}_{(\mu,a)} \times\mathbb{R}\times\mathbb{R}^{*}) $ is
given by $u_{(\mu,a)}(\Pi,\Gamma,\alpha,l)\cdot \pi_{(\mu,a)}=
u(A,c,\Pi,\Gamma, \alpha,l )|_{\mathcal{O}_{(\mu,a)}
\times\mathbb{R}\times\mathbb{R}^{*} }, $ where $\pi_{(\mu,a)}:
\mathbf{J}_Q^{-1}(\mu,a) \rightarrow \mathcal{O}_{(\mu,a)}
\times\mathbb{R}\times\mathbb{R}^{*}, \; W_{(\mu,a)}= \pi_{(\mu,a)}(W\cap
\mathbf{J}^{-1}_Q(\mu, a )). $
The $R_p$-reduced controlled spacecraft-rotor
system is the 4-tuple $(\mathcal{O}_{(\mu,a)} \times \mathbb{R}\times
\mathbb{R}^{*},\tilde{\omega}_{\mathcal{O}_{(\mu,a)} \times \mathbb{R}
\times \mathbb{R}^{*}}^{-},h_{(\mu,a)},u_{(\mu,a)}), $ where
$\tilde{\omega}_{\mathcal{O}_{(\mu,a)} \times \mathbb{R}\times
\mathbb{R}^{*}}^{-}$ is the induced symplectic form from the Poisson bracket
on $\mathfrak{se}^\ast(3) \times \mathbb{R}\times \mathbb{R}^{*} ,$ such that
Hamiltonian vector field
$$X_{h_{(\mu,a)}}(K_{(\mu,a)})=\tilde{\omega}_{\mathcal{O}_{(\mu,a)} \times \mathbb{R}\times
\mathbb{R}^{*}}^{-}(X_{K_{(\mu,a)}}, X_{h_{(\mu,a)}})
=\{K_{(\mu,a)},h_{(\mu,a)}\}_{-}|_{\mathcal{O}_{(\mu,a)}
\times\mathbb{R}\times\mathbb{R}^{*} }, $$
since $(\mathcal{O}_{(\mu,a)} \times
\mathbb{R}\times \mathbb{R}^*,\tilde{\omega}_{\mathcal{O}_{(\mu,a)}
\times \mathbb{R} \times \mathbb{R}^*}^{-})$ is a symplectic leaf of
the Poisson manifold
$(\mathfrak{se}^\ast(3) \times \mathbb{R} \times \mathbb{R}^*, \{\cdot,\cdot\}_{-}). $
Moreover, assume that the dynamical vector field of the $R_p$-reduced controlled
spacecraft-rotor system $(\mathcal{O}_{(\mu,a)} \times \mathbb{R}
\times \mathbb{R}^*,\tilde{\omega}_{\mathcal{O}_{(\mu,a)} \times
\mathbb{R} \times \mathbb{R}^*}^{-},h_{(\mu,a)},u_{(\mu,a)})$ can be
expressed by
\begin{align} X_{(\mathcal{O}_{(\mu,a)} \times \mathbb{R} \times
\mathbb{R}^*,\tilde{\omega}_{\mathcal{O}_{(\mu,a)} \times \mathbb{R}
\times \mathbb{R}^*}^{-},h_{(\mu,a)},u_{(\mu,a)})} = X_{h_{(\mu,a)}} +
\mbox{vlift}(u_{(\mu,a)}),
\label{3.11} \end{align}
where $\mbox{vlift}(u_{(\mu,a)})=
\mbox{vlift}(u_{(\mu,a)})X_{h_{(\mu,a)}} \in T(\mathcal{O}_{(\mu,a)}
\times\mathbb{R}\times\mathbb{R}^*), $ is the change of $X_{h_{(\mu,a)}}$
under the action of the $R_p$-reduced control torque $u_{(\mu,a)}$.
The dynamical vector fields of the controlled
spacecraft-rotor system and the $R_p$-reduced controlled
spacecraft-rotor system satisfy the condition
\begin{equation}X_{(\mathcal{O}_{(\mu,a)} \times \mathbb{R} \times
\mathbb{R}^*,\tilde{\omega}_{\mathcal{O}_{(\mu,a)} \times \mathbb{R} \times
\mathbb{R}^*}^{-}, h_{(\mu,a)}, u_{(\mu,a)})}\cdot \pi_{(\mu,a)}
=T\pi_{(\mu,a)}\cdot X_{(T^\ast Q,\textmd{SE}(3),\omega_Q,H,u)}\cdot i_{(\mu,a)}.
\label{3.12}\end{equation}
See Marsden et al \cite{mawazh10} and Wang \cite{wa18}.\\

To sum up the above discussion, we have the following theorem.
\begin{theo}
In the case of non-coincident centers of buoyancy and gravity, the rigid
spacecraft-rotor system with the control torque $u$ acting on the
rotor, that is, the 5-tuple $(T^\ast Q,\textmd{SE}(3),\omega_Q,H,u
), $ where $Q= \textmd{SE}(3)\times \mathbb{R}, $ is a regular point
reducible RCH system. For a point $(\mu,a) \in
\mathfrak{se}^\ast(3)$, the regular value of the momentum map
$\mathbf{J}_Q: T^* Q \cong \textmd{SE}(3)\times \mathfrak{se}^\ast(3) \times
\mathbb{R} \times \mathbb{R}^* \to \mathfrak{se}^\ast(3)$,
the $R_p$-reduced controlled spacecraft-rotor system is the 4-tuple $(\mathcal{O}_{(\mu,a)} \times
\mathbb{R} \times \mathbb{R}^*,\tilde{\omega}_{\mathcal{O}_{(\mu,a)}
\times \mathbb{R} \times \mathbb{R}^*}^{-},h_{(\mu,a)},u_{(\mu,a)}), $
where $\mathcal{O}_{(\mu,a)} \subset \mathfrak{se}^\ast(3)$ is the
co-adjoint orbit, $\tilde{\omega}_{\mathcal{O}_{(\mu,a)} \times
\mathbb{R} \times \mathbb{R}^*}^{-}$ is the induced symplectic form on
$\mathcal{O}_{(\mu,a)} \times \mathbb{R}\times \mathbb{R}^* $,
$h_{(\mu,a)}(\Pi,\Gamma,\alpha,l)\cdot \pi_{(\mu,a)}=H(A,c,\Pi,\Gamma,\alpha,l)|_{\mathcal{O}_{(\mu,a)}
\times\mathbb{R}\times\mathbb{R}^*}$, and
$u_{(\mu,a)}(\Pi,\Gamma,\alpha,l)\cdot \pi_{(\mu,a)}= u(A,c,\Pi,
\Gamma,\alpha,l)|_{\mathcal{O}_{(\mu,a)}
\times\mathbb{R}\times\mathbb{R}^*}$, and the dynamical vector field of the $R_p$-reduced controlled
spacecraft-rotor system $(\mathcal{O}_{(\mu,a)} \times \mathbb{R}
\times \mathbb{R}^*,\tilde{\omega}_{\mathcal{O}_{(\mu,a)} \times
\mathbb{R} \times \mathbb{R}^*}^{-},h_{(\mu,a)},u_{(\mu,a)})$ satisfies
(3.11) and (3.12).
\end{theo}

\section{Hamilton-Jacobi Equation of the Reduced Spacecraft-Rotor System}

The Hamilton-Jacobi theory for the regular (controlled) Hamiltonian system
is a very important subject, following the theoretical and applied development of
geometric mechanics, a lot of important problems about this subject
are being explored and studied, see Wang \cite{wa17}, Wang \cite{wa13d}
Wang \cite{wa13e}and de Le\'{o}n and Wang \cite{lewa15}.
As an application of the theoretical result,
in this section, we first give precisely the geometric constraint conditions
of the canonical symplectic form for the dynamical vector field of the controlled
rigid spacecraft-rotor system, that is, Type I and Type II of
Hamilton-Jacobi equation for the controlled rigid
spacecraft-rotor system. Then,
for the above $R_p$- reduced controlled rigid spacecraft-rotor
systems with coincident and non-coincident centers of buoyancy and gravity,
we shall derive precisely the geometric constraint conditions
of the $R_p$-reduced symplectic forms for the
dynamical vector fields of the regular point reducible controlled
rigid spacecraft-rotor systems, respectively,
that is, Type I and Type II of
Hamilton-Jacobi equations for the $R_p$-reduced controlled rigid
spacecraft-rotor systems. We shall follow the notations
and conventions introduced in Marsden \cite{ma92}, Marsden et al \cite{mawazh10},
Wang \cite{wa17} and Wang \cite{wa13d}.\\

Let $G=\textmd{SO}(3)$ or $\textmd{SE}(3)$, and $Q=\textmd{SO}(3)\times \mathbb{R}$ or
$\textmd{SE}(3)\times \mathbb{R}$, and $\omega_Q$ is canonical
symplectic form on $T^* Q$. Denote by $\Omega^i(Q)$ the set of all i-forms on $Q$, $i=1,2.$
For any $\gamma \in \Omega^1(Q),\; q\in Q, $ then $\gamma(q)\in T_q^*Q, $
and we can define a map $\gamma: Q \rightarrow T^*Q, \; q \rightarrow (q, \gamma(q)).$
Hence we say often that the map $\gamma: Q
\rightarrow T^*Q$ is an one-form on $Q$. If the one-form $\gamma$ is closed,
then $\mathbf{d}\gamma(x,y)=0, \; \forall\;
x, y \in TQ$; and the one-form $\gamma$ is called to be closed with respect to $T\pi_{Q}:
TT^* Q \rightarrow TQ, $ if for any $v, w \in TT^* Q, $ we have
$\mathbf{d}\gamma(T\pi_{Q}(v),T\pi_{Q}(w))=0. $ Since
the rigid spacecraft-rotor system with the control torque $u$ acting on the rotor
is a regular point reducible RCH system, from Theorem 2.6 and Theorem 2.7
in Wang \cite{wa13d}, we can obtain directly
the following Theorem 4.1. For convenience, the maps involved in
the following theorem are shown in Diagram-1.

\begin{center}
\hskip 0cm \xymatrix{ T^* Q \ar[r]^\varepsilon & T^* Q
\ar[d]_{X_{H\cdot \varepsilon}} \ar[dr]^{\tilde{X}^\varepsilon} \ar[r]^{\pi_Q}
& Q \ar[d]^{\tilde{X}^\gamma} \ar[r]^{\gamma}
& T^*Q \ar[d]_{\tilde{X}} \\
& T(T^*Q)  & TQ \ar[l]^{T\gamma}
& T(T^*Q) \ar[l]^{T\pi_Q} }
\end{center}
$$\mbox{Diagram-1}$$

\begin{theo}
For the controlled rigid spacecraft-rotor system $(T^*Q,\omega_Q,H,u)$ with the
canonical symplectic form $\omega_Q$ on $T^*Q$, assume that $\gamma: Q
\rightarrow T^*Q$ is an one-form on $Q$, and
$\lambda=\gamma\cdot\pi_{Q}: T^* Q \rightarrow T^* Q $,
and the map $\varepsilon: T^* Q \rightarrow T^* Q $ is symplectic.
Denote by $\tilde{X}^\gamma = T\pi_{Q}\cdot \tilde{X} \cdot \gamma$,
and $\tilde{X}^\varepsilon = T\pi_{Q}\cdot \tilde{X} \cdot \varepsilon$,
where $\tilde{X}=X_{(T^\ast Q,\omega_Q,H,u)}$ is the dynamical vector field of
the controlled rigid spacecraft-rotor system $(T^*Q,\omega_Q,H,u)$.
Then the following two assertions hold:\\
\noindent $(\mathbf{i})$
If the one-form $\gamma: Q \rightarrow T^*Q $ is closed with respect to
$T\pi_Q: TT^* Q \rightarrow TQ, $ then $\gamma$ is a solution of the Type I of
Hamilton-Jacobi equation
$T\gamma\cdot \tilde{X}^\gamma= X_H\cdot \gamma ,$ where $X_H$ is the Hamiltonian vector field
of the corresponding Hamiltonian system $(T^*Q,\omega_Q,H).$ \\
\noindent $(\mathbf{ii})$
The $\varepsilon$ is a solution of the Type II of
Hamilton-Jacobi equation $T\gamma\cdot \tilde{X}^\varepsilon= X_H\cdot
\varepsilon ,$ if and only if it is a solution of the equation
$T\varepsilon\cdot X_{H\cdot\varepsilon}= T\lambda \cdot \tilde{X} \cdot \varepsilon,$
where $X_H$ and $ X_{H\cdot\varepsilon} \in
TT^*Q $ are the Hamiltonian vector fields of the functions $H$ and $H\cdot\varepsilon:
T^*Q\rightarrow \mathbb{R}, $ respectively. \hskip 0.3cm $\blacksquare$
\end{theo}

\subsection{In The Case of Coincident Centers}

In the following we shall derive precisely the geometric constraint conditions
of the $R_p$-reduced symplectic form $\tilde{\omega}_{\mathcal{O}_\mu \times \mathbb{R} \times
\mathbb{R}^{*}}^{-}$ for the
dynamical vector field of the regular point reducible controlled
rigid spacecraft-rotor system with
coincident centers of buoyancy and gravity,
that is, Type I and Type II of
Hamilton-Jacobi equation for the $R_p$-reduced controlled rigid spacecraft-rotor system
$(\mathcal{O}_\mu \times \mathbb{R}\times \mathbb{R}^{*},
\omega^{-}_{\mathcal{O}_\mu \times \mathbb{R}\times \mathbb{R}^{*}},
h_\mu, u_\mu ) .$\\

Assume that $\gamma: \textmd{SO}(3)\times \mathbb{R} \rightarrow
T^*(\textmd{SO}(3)\times \mathbb{R})$ is an one-form on
$\textmd{SO}(3)\times \mathbb{R}$,
$\gamma(A, \alpha)=(\gamma_1, \cdots, \gamma_{8})$,
and $\gamma$ is
closed with respect to $T\pi_{\textmd{SO}(3)\times \mathbb{R}}:
TT^* (\textmd{SO}(3)\times \mathbb{R}) \rightarrow
T(\textmd{SO}(3)\times \mathbb{R}). $ For $\mu \in \mathfrak{so}^\ast(3)$
the regular value of $\mathbf{J}_Q$,
$\textmd{Im}(\gamma)\subset \mathbf{J}_Q^{-1}(\mu), $ and it is
$\textmd{SO}(3)_\mu$-invariant, and $\bar{\gamma}=\pi_\mu(\gamma):
\textmd{SO}(3)\times \mathbb{R} \rightarrow \mathcal{O}_\mu \times
\mathbb{R} \times \mathbb{R}^{*}$. Denote by
$\bar{\gamma}= (\bar{\gamma}_1,\bar{\gamma}_2,\bar{\gamma}_3,
\bar{\gamma}_4,\bar{\gamma}_5) \in \mathcal{O}_\mu \times
\mathbb{R} \times \mathbb{R}^{*}(\subset \mathfrak{so}^\ast(3)
\times \mathbb{R} \times \mathbb{R}^{*}), $ where
$\pi_\mu: \mathbf{J}_Q^{-1}(\mu) \rightarrow \mathcal{O}_\mu \times
\mathbb{R} \times \mathbb{R}^{*}. $ We choose that
$(\Pi,\alpha,l)\in
\mathcal{O}_\mu\times \mathbb{R} \times \mathbb{R}^{*}, $ and
$\Pi=(\Pi_1,\Pi_2,\Pi_3)=(\bar{\gamma}_1,\bar{\gamma}_2,\bar{\gamma}_3),
\; \alpha= \bar{\gamma}_4, \; l= \bar{\gamma}_5, $ then $h_{\mu}
\cdot \bar{\gamma}: \textmd{SO}(3)\times \mathbb{R} \rightarrow
\mathbb{R} $ is given by
\begin{align*} h_{\mu} (\Pi,\alpha,l) \cdot \bar{\gamma}  =
H (A,\Pi,\alpha,l) |_{\mathcal{O}_\mu \times
\mathbb{R} \times \mathbb{R}^{*}}\cdot \bar{\gamma}
= \frac{1}{2}[ \frac{\bar{\gamma}_1^2}{\bar{I}_1}+ \frac{\bar{\gamma}_2^2}{\bar{I}_2} + \frac{(\bar{\gamma}_3-
\bar{\gamma}_5)^2}{\bar{I}_3}+ \frac{\bar{\gamma}_5^2}{J_3}],
\end{align*} and the vector field
\begin{align*}
 X_{h_{\mu}}(\Pi) \cdot \bar{\gamma} & =\{\Pi,h_{\mu}
\}_{-}|_{\mathcal{O}_\mu \times \mathbb{R} \times
\mathbb{R}^{*}}\cdot \bar{\gamma}\\
& = -\Pi\cdot(\nabla_\Pi\Pi\times\nabla_\Pi (h_{\mu})) \cdot
\bar{\gamma}+ \{\Pi,h_{\mu} \}_{\mathbb{R}}|_{\mathcal{O}_\mu
\times \mathbb{R} \times \mathbb{R}^{*}}\cdot \bar{\gamma}\\
& = -\nabla_\Pi\Pi\cdot(\nabla_\Pi (h_{\mu})\times \Pi)\cdot
\bar{\gamma} + (\frac{\partial \Pi}{\partial \alpha}
\frac{\partial (h_\mu)}{\partial l}- \frac{\partial
(h_\mu)}{\partial \alpha}\frac{\partial \Pi}{\partial l})\cdot \bar{\gamma} \\
& =(\Pi_1,\Pi_2,\Pi_3)\times (\frac{ \Pi_1}{ \bar{I}_1},\;\;
\frac{ \Pi_2}{ \bar{I}_2}, \;\; \frac{(\Pi_3- l)}{\bar{I}_3})\cdot \bar{\gamma}\\
&= ( \frac{(\bar{I}_2-\bar{I}_3)\bar{\gamma}_2\bar{\gamma}_3-
\bar{I}_2\bar{\gamma}_2\bar{\gamma}_5 }{\bar{I}_2\bar{I}_3}, \;\;
\frac{(\bar{I}_3-\bar{I}_1)\bar{\gamma}_3\bar{\gamma}_1 +
\bar{I}_1\bar{\gamma}_1\bar{\gamma}_5}{\bar{I}_3\bar{I}_1}, \;\;
\frac{(\bar{I}_1-\bar{I}_2)\bar{\gamma}_1\bar{\gamma}_2}{\bar{I}_1\bar{I}_2} ),
\end{align*}
since $\nabla_{\Pi_i}\Pi_i=1,\; \nabla_{\Pi_i}\Pi_j=0, \; i\neq j ,$ and $\nabla_{\Pi_k} (h_{\mu})=
\Pi_k /\bar{I}_k , \; \nabla_{\Pi_3} (h_{\mu})=
(\Pi_3-l)/\bar{I}_3, $ and $\frac{\partial
\Pi}{\partial \alpha}= \frac{\partial (h_\mu)}{\partial
\alpha}=0, \; i, j=1,2,3, \; k=1,2. $

\begin{align*}
 X_{h_{\mu}}(\alpha) \cdot \bar{\gamma} & =\{\alpha,h_{\mu}
\}_{-}|_{\mathcal{O}_\mu \times \mathbb{R} \times
\mathbb{R}^{*}}\cdot \bar{\gamma}\\
& =-\Pi\cdot(\nabla_\Pi\alpha \times\nabla_\Pi (h_{\mu})) \cdot
\bar{\gamma}+ \{\alpha,h_{\mu}\}_{\mathbb{R}}|_{\mathcal{O}_\mu
\times \mathbb{R} \times \mathbb{R}^{*}}\cdot \bar{\gamma}\\
& =-\nabla_\Pi\alpha \cdot(\nabla_\Pi (h_{\mu})\times \Pi)\cdot
\bar{\gamma} + (\frac{\partial \alpha}{\partial
\alpha} \frac{\partial (h_\mu)}{\partial l}- \frac{\partial
(h_\mu)}{\partial
\alpha}\frac{\partial \alpha}{\partial l})\cdot \bar{\gamma} =
-\frac{(\bar{\gamma}_3- \bar{\gamma}_5)}{\bar{I}_3}
+\frac{\bar{\gamma}_5}{J_3},
\end{align*}
since $\nabla_{\Pi_i}\alpha=0, \; \frac{\partial \alpha}{\partial
\alpha}= 1, \; \frac{\partial (h_\mu)}{\partial \alpha}=0,$
and $\frac{\partial (h_\mu)}{\partial l}= -(\Pi_3-l)/\bar{I}_3
+\frac{l}{J_3}, \; i= 1,2,3 , $

\begin{align*}
 X_{h_{\mu}}(l) \cdot \bar{\gamma}
& =\{l,h_{\mu} \}_{-}|_{\mathcal{O}_\mu \times
\mathbb{R}\times \mathbb{R}^{*}}\cdot \bar{\gamma}\\
& = -\Pi\cdot(\nabla_\Pi l \times\nabla_\Pi (h_{\mu})) \cdot
\bar{\gamma}+ \{l, h_{\mu} \}_{\mathbb{R}}|_{\mathcal{O}_\mu
\times \mathbb{R} \times \mathbb{R}^{*}}\cdot \bar{\gamma}\\
& = -\nabla_\Pi l \cdot(\nabla_\Pi (h_{\mu})\times \Pi)\cdot
\bar{\gamma} + (\frac{\partial l}{\partial \alpha}
\frac{\partial (h_\mu)}{\partial l}- \frac{\partial
(h_\mu)}{\partial \alpha}\frac{\partial l}{\partial l})\cdot
\bar{\gamma}=0,
\end{align*}
since $\nabla_{\Pi_i} l=0,$ and $\frac{\partial l}{\partial \alpha}=
\frac{\partial (h_\mu)}{\partial \alpha}=0, \; i=1,2,3. $ \\

On the other hand, from the expressions of the dynamical vector field $\tilde{X}$
and Hamiltonian vector field $X_H$, we have that
\begin{align*}
\tilde{X}(\Pi, \alpha, l)^\gamma & =T\pi_{\textmd{SO}(3)\times \mathbb{R}}\cdot \tilde{X}\cdot\gamma(\Pi, \alpha, l)\\
& =T\pi_{\textmd{SO}(3)\times \mathbb{R}}\cdot (X_H+ \textnormal{vlift}(u))\cdot\gamma (\Pi, \alpha, l)\\
&=T\pi_{\textmd{SO}(3)\times \mathbb{R}}\cdot X_H \cdot\gamma (\Pi, \alpha, l) = X_H\cdot\gamma(\Pi, \alpha, l),
\end{align*}
that is,
\begin{align*}
 \tilde{X}(\Pi)^\gamma & = X_H(\Pi)\cdot\gamma\\
& = ( \frac{(\bar{I}_2-\bar{I}_3)\gamma_5\gamma_6-
\bar{I}_2\gamma_5\gamma_{8}}{\bar{I}_2\bar{I}_3}, \;\;
\frac{(\bar{I}_3-\bar{I}_1)\gamma_6\gamma_4+
\bar{I}_1\gamma_4\gamma_{8}}{\bar{I}_3\bar{I}_1}, \;\;
\frac{(\bar{I}_1-\bar{I}_2)\gamma_4\gamma_5}{\bar{I}_1\bar{I}_2} ),
\end{align*}
\begin{align*}
\tilde{X}( \alpha )^\gamma = X_H(\alpha)\cdot\gamma
 = -\frac{(\gamma_6- \gamma_{8})}{\bar{I}_3}
+\frac{\gamma_{8}}{J_3}, \;\;\;\;\;\;\;\;\;\;
\tilde{X}( l )^\gamma = X_H(l)\cdot\gamma =0,
\end{align*}
Since $\gamma$ is closed with respect to
$T\pi_{\textmd{SO}(3)\times \mathbb{R}}: TT^* (\textmd{SO}(3)\times \mathbb{R})
\rightarrow T(\textmd{SO}(3)\times \mathbb{R}), $
then $\pi_{\textmd{SO}(3)\times \mathbb{R}}^*(\mathbf{d}\gamma)=0.$ We choose that
$(\gamma_4,\gamma_5,\gamma_6)=\Pi=(\Pi_1,\Pi_2,\Pi_3)=
(\bar{\gamma}_1,\bar{\gamma}_2,\bar{\gamma}_3), $
and $\gamma_7= \alpha= \bar{\gamma}_4,
\; \gamma_{8}= l=\bar{\gamma}_5. $
Hence
\begin{align*}
T\bar{\gamma}\cdot \tilde{X}(\Pi)^\gamma
= X_{h_\mu}(\Pi) \cdot \bar{\gamma}, \;\;\;\;\;\;
T\bar{\gamma}\cdot \tilde{X}(\alpha)^\gamma
= X_{h_\mu}(\alpha) \cdot \bar{\gamma}, \;\;\;\;\;\;
T\bar{\gamma}\cdot \tilde{X}(l)^\gamma
= X_{h_\mu}(l) \cdot \bar{\gamma}.
\end{align*}
Thus, the Type I of Hamilton-Jacobi equation for the
$R_p$-reduced controlled rigid spacecraft-rotor system
$(\mathcal{O}_\mu\times\mathbb{R}\times\mathbb{R}^{*},
\omega_{\mathcal{O}_\mu\times\mathbb{R}\times\mathbb{R}^{*}}^{-},h_\mu, u_\mu)$ holds.\\

Next, for $\mu \in \mathfrak{so}^\ast(3),$
the regular value of $\mathbf{J}_Q$, and
a $\textmd{SO}(3)_\mu$-invariant symplectic map $\varepsilon: T^* (\textmd{SO}(3)\times \mathbb{R})
\rightarrow T^* (\textmd{SO}(3)\times \mathbb{R}),$ assume that $\varepsilon(A,\Pi, \alpha, l)
=(\varepsilon_1,\cdots, \varepsilon_{8}),$
and $\varepsilon(\mathbf{J}_Q^{-1}(\mu))\subset \mathbf{J}_Q^{-1}(\mu). $
Denote by $\bar{\varepsilon}=\pi_\mu(\varepsilon): \mathbf{J}_Q^{-1}(\mu)\rightarrow
\mathcal{O}_\mu\times \mathbb{R} \times \mathbb{R}^{*}, $ and
$\bar{\varepsilon}=(\bar{\varepsilon}_1,\bar{\varepsilon}_2,\bar{\varepsilon}_3, \bar{\varepsilon}_4, \bar{\varepsilon}_5) \in
\mathcal{O}_\mu\times \mathbb{R}\times \mathbb{R}^{*} (\subset \mathfrak{so}^\ast(3)
\times \mathbb{R}\times \mathbb{R}^{*}), $ and
$\lambda= \gamma \cdot \pi_{\textmd{SO}(3)\times \mathbb{R}}: T^* (\textmd{SO}(3)\times \mathbb{R})
\rightarrow T^* (\textmd{SO}(3)\times \mathbb{R}),$ and $\lambda(A,\Pi, \alpha, l)
=(\lambda_1,\cdots, \lambda_{8}),$
and $\bar{\lambda}=\pi_\mu(\lambda): \mathbf{J}_Q^{-1}(\mu) \rightarrow
\mathcal{O}_\mu\times \mathbb{R}\times \mathbb{R}^{*}, $ and
$\bar{\lambda}= (\bar{\lambda}_1,\bar{\lambda}_2,\bar{\lambda}_3, \bar{\lambda}_4, \bar{\lambda}_5 ) \in
\mathcal{O}_\mu\times \mathbb{R} \times \mathbb{R}^{*}. $ We choose that
$(\Pi,\alpha,l)\in
\mathcal{O}_\mu\times \mathbb{R} \times \mathbb{R}^{*}, $ and
$\Pi=(\Pi_1,\Pi_2,\Pi_3)=(\bar{\varepsilon}_1,\bar{\varepsilon}_2,\bar{\varepsilon}_3),$
and $\alpha= \bar{\varepsilon}_4,$ and $ l=\bar{\varepsilon}_5,$
then $h_{\mu} \cdot \bar{\varepsilon}: T^*(\textmd{SO}(3)\times \mathbb{R})
\rightarrow \mathbb{R} $ is given by
\begin{align*} h_{\mu}(\Pi,\alpha,l) \cdot \bar{\varepsilon}=
H(A,\Pi, \alpha, l)|_{\mathcal{O}_\mu\times \mathbb{R} \times \mathbb{R}^{*}} \cdot \bar{\varepsilon}
= \frac{1}{2}[ \frac{\bar{\varepsilon}_1^2}{\bar{I}_1}+ \frac{\bar{\varepsilon}_2^2}{\bar{I}_2}
+\frac{(\bar{\varepsilon}_3-
\bar{\varepsilon}_5)^2}{\bar{I}_3}+ \frac{\bar{\varepsilon}_5^2}{J_3}],
\end{align*}
and the vector field
\begin{align*}
 X_{h_{\mu}}(\Pi) \cdot \bar{\varepsilon}
& =\{\Pi,h_{\mu} \}_{-}|_{\mathcal{O}_\mu\times \mathbb{R} \times \mathbb{R}^{*}} \cdot
\bar{\varepsilon}\\
& = -\Pi\cdot(\nabla_\Pi\Pi\times\nabla_\Pi (h_{\mu})) \cdot
\bar{\varepsilon}+ \{\Pi,h_{\mu} \}_{\mathbb{R}}|_{\mathcal{O}_\mu
\times \mathbb{R} \times \mathbb{R}^{*}}\cdot \bar{\varepsilon}\\
&= ( \frac{(\bar{I}_2-\bar{I}_3)\bar{\varepsilon}_2\bar{\varepsilon}_3-
\bar{I}_2\bar{\varepsilon}_2\bar{\varepsilon}_5 }{\bar{I}_2\bar{I}_3}, \;\;
\frac{(\bar{I}_3-\bar{I}_1)\bar{\varepsilon}_3\bar{\varepsilon}_1 +
\bar{I}_1\bar{\varepsilon}_1\bar{\varepsilon}_5}{\bar{I}_3\bar{I}_1}, \;\;
\frac{(\bar{I}_1-\bar{I}_2)\bar{\varepsilon}_1\bar{\varepsilon}_2}{\bar{I}_1\bar{I}_2} ),
\end{align*}
\begin{align*}
 X_{h_{\mu}}(\alpha) \cdot \bar{\varepsilon} & =\{\alpha,h_{\mu}
\}_{-}|_{\mathcal{O}_\mu \times \mathbb{R} \times
\mathbb{R}^{*}}\cdot \bar{\varepsilon}\\
& =-\Pi\cdot(\nabla_\Pi\alpha \times\nabla_\Pi (h_{\mu})) \cdot
\bar{\varepsilon}+ \{\alpha,h_{\mu}\}_{\mathbb{R}}|_{\mathcal{O}_\mu
\times \mathbb{R} \times \mathbb{R}^{*}}\cdot \bar{\varepsilon}
 = -\frac{(\bar{\varepsilon}_3- \bar{\varepsilon}_5)}{\bar{I}_3}
+\frac{\bar{\varepsilon}_5}{J_3},
\end{align*}
\begin{align*}
 X_{h_{\mu}}(l) \cdot \bar{\varepsilon}
& =\{l,h_{\mu} \}_{-}|_{\mathcal{O}_\mu \times
\mathbb{R} \times \mathbb{R}^{*}}\cdot \bar{\varepsilon}\\
& = -\Pi\cdot(\nabla_\Pi l \times\nabla_\Pi (h_{\mu})) \cdot
\bar{\varepsilon}+ \{l, h_{\mu} \}_{\mathbb{R}}|_{\mathcal{O}_\mu
\times \mathbb{R} \times \mathbb{R}^{*}}\cdot \bar{\varepsilon}=0.
\end{align*}

On the other hand, from the expressions of the dynamical vector field $\tilde{X}$
and Hamiltonian vector field $X_H$, we have that
\begin{align*}
\tilde{X}(\Pi, \alpha, l)^\varepsilon & =T\pi_{\textmd{SO}(3)\times \mathbb{R}}\cdot \tilde{X}\cdot\varepsilon(\Pi, \alpha, l)\\
& =T\pi_{\textmd{SO}(3)\times \mathbb{R}}\cdot (X_H+ \textnormal{vlift}(u))\cdot\varepsilon (\Pi, \alpha, l)\\
& =T\pi_{\textmd{SO}(3)\times \mathbb{R}}\cdot X_H \cdot\varepsilon (\Pi, \alpha, l)= X_H\cdot\varepsilon(\Pi, \alpha, l),
\end{align*}
that is,
\begin{align*}
 \tilde{X}(\Pi)^\varepsilon  & = X_H(\Pi)\cdot\varepsilon\\
& = ( \frac{(\bar{I}_2-\bar{I}_3)\varepsilon_5\varepsilon_6-
\bar{I}_2\varepsilon_5\varepsilon_{8}}{\bar{I}_2\bar{I}_3}, \;\;
\frac{(\bar{I}_3-\bar{I}_1)\varepsilon_6\varepsilon_4+
\bar{I}_1\varepsilon_4\varepsilon_{8}}{\bar{I}_3\bar{I}_1}, \;\;
\frac{(\bar{I}_1-\bar{I}_2)\varepsilon_4\varepsilon_5}{\bar{I}_1\bar{I}_2} ),
\end{align*}
\begin{align*}
\tilde{X}( \alpha )^\varepsilon = X_H(\alpha)\cdot\varepsilon
 = -\frac{(\varepsilon_6- \varepsilon_{8})}{\bar{I}_3}
+\frac{\varepsilon_{8}}{J_3}, \;\;\;\;\;\;\;\;\;\;
\tilde{X}( l )^\varepsilon = X_H(l)\cdot\varepsilon =0.
\end{align*}
Note that
\begin{align*}
T\bar{\gamma}\cdot \tilde{X}(\Pi)^\varepsilon
= ( \frac{(\bar{I}_2-\bar{I}_3)\bar{\gamma}_2\bar{\gamma}_3-
\bar{I}_2\bar{\gamma}_2\bar{\gamma}_5}{\bar{I}_2\bar{I}_3}, \;\;
\frac{(\bar{I}_3-\bar{I}_1)\bar{\gamma}_3\bar{\gamma}_1+
\bar{I}_1\bar{\gamma}_1\bar{\gamma}_5}{\bar{I}_3\bar{I}_1}, \;\;
\frac{(\bar{I}_1-\bar{I}_2)\bar{\gamma}_1\bar{\gamma}_2}{\bar{I}_1\bar{I}_2} ),
\end{align*}
\begin{align*}
T\bar{\gamma}\cdot \tilde{X}(\alpha)^\varepsilon
& = -\frac{(\bar{\gamma}_3- \bar{\gamma}_5)}{\bar{I}_3}
+\frac{\bar{\gamma}_5}{J_3}, \;\;\;\;\;\;\;\;\;\;
T\bar{\gamma}\cdot \tilde{X}(l)^\varepsilon= 0,
\end{align*}
and $$T\bar{\lambda}\cdot \tilde{X} \cdot \varepsilon=T\pi_\mu\cdot T\lambda \cdot (X_H+ \textnormal{vlift}(u))\cdot\varepsilon
=T\pi_\mu\cdot T\gamma \cdot T\pi_{\textmd{SO}(3)\times \mathbb{R}}\cdot (X_H+ \textnormal{vlift}(u))\cdot\varepsilon
=T\bar{\lambda}\cdot X_H \cdot \varepsilon,$$ that is,
\begin{align*}
 T\bar{\lambda}\cdot \tilde{X}(\Pi) \cdot \varepsilon &=T\bar{\lambda}\cdot X_H(\Pi) \cdot \varepsilon \\
&= ( \frac{(\bar{I}_2-\bar{I}_3)\bar{\lambda}_2\bar{\lambda}_3-
\bar{I}_2\bar{\lambda}_2\bar{\lambda}_5}{\bar{I}_2\bar{I}_3}, \;\;
\frac{(\bar{I}_3-\bar{I}_1)\bar{\lambda}_3\bar{\lambda}_1+
\bar{I}_1\bar{\lambda}_1\bar{\lambda}_5}{\bar{I}_3\bar{I}_1}, \;\;
\frac{(\bar{I}_1-\bar{I}_2)\bar{\lambda}_1\bar{\lambda}_2}{\bar{I}_1\bar{I}_2} ),
\end{align*}
\begin{align*}
 T\bar{\lambda}\cdot \tilde{X}(\alpha) \cdot \varepsilon=T\bar{\lambda}\cdot X_H(\alpha) \cdot \varepsilon
 = -\frac{(\bar{\lambda}_3- \bar{\lambda}_5)}{\bar{I}_3}
+\frac{\bar{\lambda}_5}{J_3},
\end{align*}
\begin{align*}
T\bar{\lambda}\cdot \tilde{X}(l) \cdot \varepsilon=T\bar{\lambda}\cdot X_H(l) \cdot \varepsilon=0.
\end{align*}
Thus, when we choose that $(\Pi,\alpha,l)\in
\mathcal{O}_\mu\times \mathbb{R} \times \mathbb{R}^{*}, $ and
$\Pi=(\Pi_1,\Pi_2,\Pi_3)=(\bar{\gamma}_1,\bar{\gamma}_2,\bar{\gamma}_3)=
(\bar{\varepsilon}_1,\bar{\varepsilon}_2,\bar{\varepsilon}_3)=
(\bar{\lambda}_1,\bar{\lambda}_2,\bar{\lambda}_3), $
and $\alpha= \bar{\varepsilon}_4=\bar{\lambda}_4,$ and
$ l=\bar{\varepsilon}_5= \bar{\lambda}_5,$
we must have that
\begin{align*}
& T\bar{\gamma}\cdot \tilde{X}(\Pi)^\varepsilon=X_{h_{\mu}}(\Pi) \cdot \bar{\varepsilon}
=T\bar{\lambda}\cdot \tilde{X}(\Pi) \cdot \varepsilon, \\
& T\bar{\gamma}\cdot \tilde{X}(\alpha)^\varepsilon=X_{h_{\mu}}(\alpha) \cdot \bar{\varepsilon}
=T\bar{\lambda}\cdot \tilde{X}(\alpha) \cdot \varepsilon, \\
& T\bar{\gamma}\cdot \tilde{X}(l)^\varepsilon=X_{h_{\mu}}(l) \cdot \bar{\varepsilon}
=T\bar{\lambda}\cdot \tilde{X}(l) \cdot \varepsilon.
\end{align*}
Since the map $\varepsilon: T^* (\textmd{SO}(3)\times \mathbb{R})
\rightarrow T^* (\textmd{SO}(3)\times \mathbb{R})$ is symplectic, then
$T\bar{\varepsilon}\cdot X_{h_{\mu} \cdot \bar{\varepsilon}}
=X_{h_{\mu}} \cdot \bar{\varepsilon}. $
Thus, in this case, we must have that
$\varepsilon$ and $\bar{\varepsilon} $ are the solution of the Type II of
Hamilton-Jacobi equation
$T\bar{\gamma}\cdot \tilde{X}^\varepsilon= X_{h_{\mu}}\cdot \bar{\varepsilon}, $
for the $R_p$-reduced controlled rigid spacecraft-rotor system
$(\mathcal{O}_\mu\times \mathbb{R}\times \mathbb{R}^{*},
\omega_{\mathcal{O}_\mu\times \mathbb{R} \times \mathbb{R}^{*}}^{-},h_{\mu}, u_{\mu}) $,
if and only if they satisfy
the equation $T\bar{\varepsilon}\cdot(X_{h_{\mu} \cdot \bar{\varepsilon}})
= T\bar{\lambda}\cdot \tilde{X} \cdot\varepsilon. $\\

To sum up the above discussion, we have the following Theorem 4.2.
For convenience, the maps involved in
the following theorem are shown in Diagram-2.

\begin{center}
\hskip 0cm \xymatrix{ \mathbf{J}_Q^{-1}(\mu) \ar[r]^{i_\mu} & T^* Q
\ar[d]_{X_{H\cdot \varepsilon}} \ar[dr]^{\tilde{X}^\varepsilon} \ar[r]^{\pi_Q}
& Q \ar[d]^{\tilde{X}^\gamma} \ar[r]^{\gamma}
& T^*Q \ar[d]_{\tilde{X}} \ar[dr]_{X_{h_\mu \cdot\bar{\varepsilon}}} \ar[r]^{\pi_\mu}
& \;\; \mathcal{O}_\mu\times \mathbb{R}\times \mathbb{R}^{*} \ar[d]^{X_{h_\mu}} \\
& T(T^*Q)  & TQ \ar[l]^{T\gamma}
& T(T^*Q) \ar[l]^{T\pi_Q} \ar[r]_{T\pi_\mu}
& \;\; T(\mathcal{O}_\mu\times \mathbb{R}\times \mathbb{R}^{*})}
\end{center}
$$\mbox{Diagram-2}$$

\begin{theo}
In the case of coincident centers of buoyancy and gravity,
if the 5-tuple $(T^\ast Q,\\ \textmd{SO}(3),\omega_Q,H,u), $ where $Q=
\textmd{SO}(3)\times \mathbb{R}, $ is a regular point reducible
rigid spacecraft-rotor system with the control torque $u$ acting on the rotor,
then for a point $\mu \in \mathfrak{so}^\ast(3)$, the regular
value of the momentum map $\mathbf{J}_Q: \textmd{SO}(3)\times
\mathfrak{so}^\ast(3) \times \mathbb{R}\times \mathbb{R}^{*} \to
\mathfrak{so}^\ast(3)$, the $R_p$-reduced controlled rigid spacecraft-rotor system is the 4-tuple
$(\mathcal{O}_\mu \times \mathbb{R} \times
\mathbb{R}^{*},\tilde{\omega}_{\mathcal{O}_\mu \times \mathbb{R}
\times \mathbb{R}^{*}}^{-},h_\mu,u_\mu). $ Assume that $\gamma:
\textmd{SO}(3)\times \mathbb{R} \rightarrow
T^*(\textmd{SO}(3)\times \mathbb{R})$ is an one-form on
$\textmd{SO}(3)\times \mathbb{R}$,
and $\lambda=\gamma \cdot \pi_{(\textmd{SO}(3)\times \mathbb{R})}:
T^* (\textmd{SO}(3)\times \mathbb{R})\rightarrow
T^* (\textmd{SO}(3)\times \mathbb{R}), $ and $\varepsilon:
T^* (\textmd{SO}(3) \times \mathbb{R})\rightarrow
T^* (\textmd{SO}(3)\times \mathbb{R}) $ is a
$\textmd{SO}(3)_\mu$-invariant symplectic map.
Denote
$\tilde{X}^\gamma = T\pi_{(\textmd{SO}(3)\times \mathbb{R})}\cdot \tilde{X}\cdot \gamma$, and
$\tilde{X}^\varepsilon = T\pi_{(\textmd{SO}(3)\times \mathbb{R})}\cdot \tilde{X}\cdot \varepsilon$,
where $\tilde{X}=X_{(T^\ast Q,\textmd{SO}(3),\omega_Q,H,u)}$ is the
dynamical vector field of the controlled rigid
spacecraft-rotor system $(T^\ast Q,\textmd{SO}(3),\omega_Q,H,u)$.
Moreover, assume that $\textmd{Im}(\gamma)\subset \mathbf{J}_Q^{-1}(\mu), $ and it is
$\textmd{SO}(3)_\mu$-invariant,
and $\varepsilon(\mathbf{J}_Q^{-1}(\mu))\subset \mathbf{J}_Q^{-1}(\mu). $
Denote $\bar{\gamma}=\pi_\mu(\gamma):
\textmd{SO}(3)\times \mathbb{R}\rightarrow \mathcal{O}_\mu \times \mathbb{R}\times \mathbb{R}^{*}, $ and
$\bar{\lambda}=\pi_\mu(\lambda): T^* (\textmd{SO}(3)\times \mathbb{R}) \rightarrow
\mathcal{O}_\mu\times \mathbb{R}\times \mathbb{R}^{*}, $ and
$\bar{\varepsilon}=\pi_\mu(\varepsilon): \mathbf{J}_Q^{-1}(\mu)\rightarrow
\mathcal{O}_\mu\times \mathbb{R}\times \mathbb{R}^{*}. $
Then the following two assertions hold:\\
\noindent $(\mathbf{i})$
If the one-form $\gamma: \textmd{SO}(3)\times \mathbb{R} \rightarrow
T^*(\textmd{SO}(3)\times \mathbb{R}) $ is closed with respect to
$T\pi_{(\textmd{SO}(3)\times \mathbb{R})}: TT^* (\textmd{SO}(3)\times \mathbb{R})
\rightarrow T(\textmd{SO}(3)\times \mathbb{R}), $
then $\bar{\gamma}$ is a solution of the Type I of Hamilton-Jacobi equation
$T\bar{\gamma}\cdot \tilde{X}^\gamma= X_{h_\mu}\cdot \bar{\gamma}; $\\
\noindent $(\mathbf{ii})$
The $\varepsilon$ and $\bar{\varepsilon} $ satisfy the Type II of Hamilton-Jacobi equation
$T\bar{\gamma}\cdot \tilde{X}^\varepsilon= X_{h_\mu}\cdot \bar{\varepsilon}, $
if and only if they satisfy
the equation $T\bar{\varepsilon}\cdot(X_{h_\mu \cdot \bar{\varepsilon}})
= T\bar{\lambda}\cdot \tilde{X}\cdot\varepsilon. $ \hskip 0.3cm $\blacksquare$
\end{theo}

\begin{rema}
When the rigid spacecraft does not carry any internal rotor, in this
case the configuration space is $Q=G=\textmd{SO}(3), $
the motion of rigid spacecraft is just the rotation motion of a rigid body, the
above $R_p$-reduced controlled spacecraft-rotor system is just the
Marsden-Weinstein reduced rigid body system, that is,
3-tuple $(\mathcal{O}_{\mu},
\omega_{\mathcal{O}_{\mu}},h_{\mathcal{O}_{\mu}})$, where
$\mathcal{O}_{\mu} \subset \mathfrak{so}^\ast(3)$ is the
co-adjoint orbit, $\omega_{\mathcal{O}_{\mu}}$ is the orbit
symplectic form on $\mathcal{O}_{\mu}$, which is induced by the
rigid body Lie-Poisson bracket on $\mathfrak{so}^\ast(3)$,
$h_{\mathcal{O}_{\mu}}(\Pi)\cdot \pi_{\mathcal{O}_{\mu}}
=H(A,\Pi)|_{\mathcal{O}_{\mu}}$. From the above Theorem 4.2
we can obtain the Proposition 5.3 in Wang \cite{wa17}, that is, we give the two
types of Lie-Poisson Hamilton-Jacobi equation for the Marsden-Weinstein
reduced rigid body system $(\mathcal{O}_{\mu},
\omega_{\mathcal{O}_{\mu}},h_{\mathcal{O}_{\mu}})$.
See Marsden and Ratiu \cite{mara99},
Ge and Marsden \cite{gema88}, and Wang \cite{wa17}.
\end{rema}

It is worthy of noting that, for the controlled rigid
spacecraft-rotor system $(T^\ast Q,\textmd{SO}(3),\omega_Q,H,u)$
with the $R_p$-reduced controlled rigid spacecraft-rotor system
$(\mathcal{O}_\mu \times \mathbb{R} \times
\mathbb{R}^{*},\tilde{\omega}_{\mathcal{O}_\mu \times \mathbb{R}
\times \mathbb{R}^{*}}^{-}, h_\mu, u_\mu) $, we know that the Hamiltonian vector fields
$X_{H}$ and $X_{h_\mu}$ for the corresponding Hamiltonian system $(T^*Q,\textmd{SO}(3),\omega_Q, H)$
and its $R_p$-reduced system $(\mathcal{O}_\mu \times \mathbb{R} \times
\mathbb{R}^{*},\tilde{\omega}_{\mathcal{O}_\mu \times \mathbb{R}
\times \mathbb{R}^{*}}^{-}, h_\mu )$, are $\pi_\mu$-related, that is,
$X_{h_\mu}\cdot \pi_\mu=T\pi_\mu\cdot X_{H}\cdot i_\mu.$ Then we can
prove the following Theorem 4.4, which states the relationship
between the solutions of Type II of Hamilton-Jacobi equations and the
regular point reduction.

\begin{theo}
In the case of coincident centers of buoyancy and gravity, for the controlled rigid
spacecraft-rotor system $(T^\ast Q,\textmd{SO}(3),\omega_Q,H,u)$
with the $R_p$-reduced controlled rigid spacecraft-rotor system
$(\mathcal{O}_\mu \times \mathbb{R} \times
\mathbb{R}^{*},\tilde{\omega}_{\mathcal{O}_\mu \times \mathbb{R}
\times \mathbb{R}^{*}}^{-}, h_\mu, u_\mu) $,
assume that $\gamma:
\textmd{SO}(3)\times \mathbb{R} \rightarrow
T^*(\textmd{SO}(3)\times \mathbb{R})$ is an one-form on
$\textmd{SO}(3)\times \mathbb{R}$, and $\varepsilon:
T^* (\textmd{SO}(3) \times \mathbb{R})\rightarrow
T^* (\textmd{SO}(3)\times \mathbb{R}) $ is a
$\textmd{SO}(3)_\mu$-invariant symplectic map,
$\bar{\varepsilon}=\pi_\mu(\varepsilon): \mathbf{J}_Q^{-1}(\mu)\rightarrow
\mathcal{O}_\mu\times \mathbb{R}\times \mathbb{R}^{*}. $
Under the hypotheses and notations of Theorem 4.2, then we have that
$\varepsilon$ is a solution of the Type II of Hamilton-Jacobi equation
$T\gamma\cdot \tilde{X}^\varepsilon= X_H\cdot \varepsilon, $ for the
regular point reducible controlled rigid
spacecraft-rotor system $(T^\ast Q,\textmd{SO}(3),\omega_Q,H,u), $ if and only if
$\varepsilon$ and $\bar{\varepsilon} $ satisfy the Type II of Hamilton-Jacobi equation
$T\bar{\gamma}\cdot \tilde{X}^\varepsilon= X_{h_\mu}\cdot \bar{\varepsilon}, $ for the
$R_p$-reduced controlled rigid spacecraft-rotor system
$(\mathcal{O}_\mu \times \mathbb{R} \times
\mathbb{R}^{*},\tilde{\omega}_{\mathcal{O}_\mu \times \mathbb{R}
\times \mathbb{R}^{*}}^{-}, h_\mu, u_\mu) $.
\end{theo}

\noindent{\bf Proof: }
Note that $\textmd{Im}(\gamma)\subset \mathbf{J}_Q^{-1}(\mu), $ and it
is $\textmd{SO}(3)_\mu$-invariant, in this case,
$\pi_\mu^*\tilde{\omega}_{\mathcal{O}_\mu \times \mathbb{R}
\times \mathbb{R}^{*}}^{-}= i_\mu^*\omega_Q= \omega_Q, $ along $\textmd{Im}(\gamma)$.
Since the Hamiltonian vector fields
$X_{H}$ and $X_{h_\mu}$ are $\pi_\mu$-related, that is,
$X_{h_\mu}\cdot \pi_\mu= T\pi_\mu\cdot X_{H}\cdot i_\mu, $ and
by using the $R_p$-reduced symplectic form $\tilde{\omega}_{\mathcal{O}_\mu \times \mathbb{R}
\times \mathbb{R}^{*}}^{-}$, for any $w \in TT^* Q,$ and $T\pi_{\mu} \cdot w
\neq 0, $ we have that
\begin{align*}
& \tilde{\omega}_{\mathcal{O}_\mu \times \mathbb{R}
\times \mathbb{R}^{*}}^{-}(T\bar{\gamma} \cdot \tilde{X}^\varepsilon
- X_{h_\mu} \cdot \bar{\varepsilon}, \; T\pi_\mu \cdot w) \\
& = \tilde{\omega}_{\mathcal{O}_\mu \times \mathbb{R}
\times \mathbb{R}^{*}}^{-}(T\bar{\gamma} \cdot \tilde{X}^\varepsilon, \; T\pi_\mu \cdot w)-
\tilde{\omega}_{\mathcal{O}_\mu \times \mathbb{R}
\times \mathbb{R}^{*}}^{-}(X_{h_\mu} \cdot \bar{\varepsilon}, \; T\pi_\mu \cdot w) \\
& = \tilde{\omega}_{\mathcal{O}_\mu \times \mathbb{R}
\times \mathbb{R}^{*}}^{-}(T\pi_\mu \cdot T\gamma \cdot \tilde{X}^\varepsilon, \; T\pi_\mu \cdot w)-
\tilde{\omega}_{\mathcal{O}_\mu \times \mathbb{R}
\times \mathbb{R}^{*}}^{-}(X_{h_\mu} \cdot \pi_\mu \cdot \varepsilon, \; T\pi_\mu \cdot w) \\
& = \pi_\mu^*\tilde{\omega}_{\mathcal{O}_\mu \times \mathbb{R}
\times \mathbb{R}^{*}}^{-}(T\gamma \cdot \tilde{X}^\varepsilon, \; w)
-\tilde{\omega}_{\mathcal{O}_\mu \times \mathbb{R}
\times \mathbb{R}^{*}}^{-}(T\pi_\mu\cdot X_{H}\cdot \varepsilon, \; T\pi_\mu \cdot w) \\
& = \pi_\mu^*\tilde{\omega}_{\mathcal{O}_\mu \times \mathbb{R}
\times \mathbb{R}^{*}}^{-}(T\gamma \cdot \tilde{X}^\varepsilon, \; w)
-\pi_\mu^*\tilde{\omega}_{\mathcal{O}_\mu \times \mathbb{R}
\times \mathbb{R}^{*}}^{-}(X_{H}\cdot \varepsilon, \; w)\\
& = \omega_Q(T\gamma \cdot \tilde{X}^\varepsilon- X_{H}\cdot \varepsilon, \; w).
\end{align*}
Because both the symplectic form $\omega_Q$ and the $R_p$-reduced symplectic form
$\tilde{\omega}_{\mathcal{O}_\mu \times \mathbb{R}\times \mathbb{R}^{*}}^{-}$
are non-degenerate, it follows that the equation
$T\bar{\gamma}\cdot \tilde{X}^\varepsilon= X_{h_\mu}\cdot \bar{\varepsilon}, $
is equivalent to the equation $T\gamma\cdot \tilde{X}^\varepsilon= X_H\cdot \varepsilon$. Thus,
$\varepsilon$ is a solution of the Type II of Hamilton-Jacobi equation
$T\gamma\cdot \tilde{X}^\varepsilon= X_H\cdot \varepsilon, $ for the
regular point reducible controlled rigid
spacecraft-rotor system $(T^\ast Q,\textmd{SO}(3),\omega_Q,H,u), $ if and only if
$\varepsilon$ and $\bar{\varepsilon} $ satisfy the Type II of Hamilton-Jacobi equation
$T\bar{\gamma}\cdot \tilde{X}^\varepsilon= X_{h_\mu}\cdot \bar{\varepsilon}, $ for the
$R_p$-reduced controlled rigid spacecraft-rotor system
$(\mathcal{O}_\mu \times \mathbb{R} \times
\mathbb{R}^{*},\tilde{\omega}_{\mathcal{O}_\mu \times \mathbb{R}
\times \mathbb{R}^{*}}^{-}, h_\mu, u_\mu) $.  \hskip 0.3cm
$\blacksquare$

\subsection{In The Case of Non-coincident Centers}

In the following we shall derive precisely the geometric constraint conditions
of the $R_p$-reduced symplectic form $\tilde{\omega}_{\mathcal{O}_{(\mu,a)} \times \mathbb{R} \times
\mathbb{R}^{*}}^{-}$ for the
dynamical vector field of the regular point reducible controlled
rigid spacecraft-rotor system with
non-coincident centers of buoyancy and gravity,
that is, Type I and Type II of
Hamilton-Jacobi equation for the $R_p$-reduced controlled rigid spacecraft-rotor system
$(\mathcal{O}_{(\mu,a)}\times \mathbb{R}\times \mathbb{R}^{*},
\omega^{-}_{\mathcal{O}_{(\mu,a)}\times \mathbb{R}\times \mathbb{R}^{*}},
h_{(\mu,a)}, u_{(\mu,a)}) .$\\

Assume that $\gamma: \textmd{SE}(3)\times \mathbb{R} \rightarrow
T^* (\textmd{SE}(3) \times \mathbb{R}) $ is an one-form on
$\textmd{SE}(3) \times \mathbb{R} $, and
$\gamma(A,c,\alpha)=(\gamma_1, \cdots, \gamma_{14})$,
and $\gamma$ is closed with respect to $T\pi_{\textmd{SE}(3)\times \mathbb{R}}:
TT^* (\textmd{SE}(3) \times \mathbb{R}) \rightarrow
T(\textmd{SE}(3) \times \mathbb{R}). $
For $(\mu,a) \in \mathfrak{se}^\ast(3),$ the regular value of $\mathbf{J}_Q$,
$\textmd{Im}(\gamma)\subset \mathbf{J}_Q^{-1}(\mu,a), $ and it is
$\textmd{SE}(3)_{(\mu,a)}$-invariant, and $\bar{\gamma}=\pi_{(\mu,a)}(\gamma):
\textmd{SE}(3)\times \mathbb{R}\rightarrow \mathcal{O}_{(\mu,a)} \times
\mathbb{R}\times \mathbb{R}^{*}$. Denote by
$\bar{\gamma}= (\bar{\gamma}_1,\bar{\gamma}_2,\bar{\gamma}_3,
\bar{\gamma}_4,\bar{\gamma}_5, \bar{\gamma}_6,\bar{\gamma}_7,
\bar{\gamma}_8 ) \in \mathcal{O}_{(\mu,a)} \times
\mathbb{R}\times \mathbb{R}^{*}(\subset \mathfrak{se}^\ast(3)
\times \mathbb{R}\times \mathbb{R}^{*}), $ where
$\pi_{(\mu,a)}: \mathbf{J}_Q^{-1}(\mu,a) \rightarrow \mathcal{O}_{(\mu,a)} \times
\mathbb{R} \times \mathbb{R}^{*}. $ We choose that
$(\Pi,\Gamma,\alpha,l)\in
\mathcal{O}_{(\mu,a)}\times \mathbb{R} \times \mathbb{R}^{*}, $ and
$\Pi=(\Pi_1,\Pi_2,\Pi_3)=(\bar{\gamma}_1,\bar{\gamma}_2,\bar{\gamma}_3)$,
$\Gamma=(\Gamma_1,\Gamma_2,\Gamma_3)=(\bar{\gamma}_4,\bar{\gamma}_5,\bar{\gamma}_6)$,
$\alpha=\bar{\gamma}_7, $ and $ l= \bar{\gamma}_8. $ Then $h_{(\mu,a)} \cdot \bar{\gamma}:
\textmd{SE}(3) \times \mathbb{R} \rightarrow \mathbb{R} $ is given
by
\begin{align*} & h_{(\mu,a)}(\Pi,\Gamma,\alpha,l) \cdot \bar{\gamma}=
H (A, c, \Pi,\Gamma,\alpha,l) |_{\mathcal{O}_{(\mu,a)} \times \mathbb{R} \times \mathbb{R}^{*}}
\cdot \bar{\gamma}\\ &
=\frac{1}{2}[ \frac{\bar{\gamma}_1^2}{\bar{I}_1}+ \frac{\bar{\gamma}_2^2}{\bar{I}_2}
+\frac{(\bar{\gamma}_3-\bar{\gamma}_8)^2}{\bar{I}_3}+ \frac{\bar{\gamma}_8^2}{J_3}]
+ gh(\bar{\gamma}_4\cdot\chi_1+ \bar{\gamma}_5\cdot\chi_2+
\bar{\gamma}_6\cdot\chi_3), \end{align*} and the vector field
\begin{align*}
 X_{h_{(\mu,a)}}(\Pi) \cdot \bar{\gamma} & =\{\Pi,\;
h_{(\mu,a)}\}_{-}|_{\mathcal{O}_{(\mu,a)}
\times \mathbb{R}\times \mathbb{R}^{*}}\cdot \bar{\gamma}\\
& = -\Pi\cdot(\nabla_\Pi\Pi\times\nabla_\Pi (h_{(\mu,a)})) \cdot
\bar{\gamma}-\Gamma\cdot(\nabla_\Pi\Pi\times\nabla_\Gamma
(h_{(\mu,a)})-\nabla_\Pi
(h_{(\mu,a)}) \times\nabla_\Gamma\Pi)\cdot \bar{\gamma}\\
& \;\;\;\; + \{\Pi,\;
h_{(\mu,a)}\}_{\mathbb{R}}|_{\mathcal{O}_{(\mu,a)} \times
\mathbb{R}\times \mathbb{R}^{*}}\cdot \bar{\gamma}\\
& = -\nabla_\Pi\Pi \cdot(\nabla_\Pi (h_{(\mu,a)}) \times \Pi)\cdot
\bar{\gamma}- \nabla_\Pi\Pi \cdot (\nabla_\Gamma
(h_{(\mu,a)}) \times \Gamma)\cdot \bar{\gamma}\\
& \;\;\;\; + (\frac{\partial \Pi}{\partial \alpha}
\frac{\partial (h_{(\mu,a)})}{\partial l}- \frac{\partial
(h_{(\mu,a)})}{\partial
\alpha}\frac{\partial \Pi}{\partial l})\cdot \bar{\gamma}\\
&=(\Pi_1,\Pi_2,\Pi_3)\times (\frac{ \Pi_1}{ \bar{I}_1},
\frac{ \Pi_2}{ \bar{I}_2}, \frac{(\Pi_3-l)}{\bar{I}_3})\cdot
\bar{\gamma}
+gh(\Gamma_1,\Gamma_2,\Gamma_3)\times (\chi_1,\chi_2,\chi_3)\cdot \bar{\gamma}\\
&= ( \frac{(\bar{I}_2-\bar{I}_3)\bar{\gamma}_2\bar{\gamma}_3-
\bar{I}_2\bar{\gamma}_2\bar{\gamma}_{8}}{\bar{I}_2\bar{I}_3}+
gh(\bar{\gamma}_5\chi_3-\bar{\gamma}_6\chi_2), \\
& \;\;\;\;\;\; \frac{(\bar{I}_3-\bar{I}_1)
\bar{\gamma}_3 \bar{\gamma}_1+
\bar{I}_1\bar{\gamma}_1\bar{\gamma}_8}{\bar{I}_3\bar{I}_1} +
gh(\bar{\gamma}_6\chi_1-\bar{\gamma}_4\chi_3), \;\;
\frac{(\bar{I}_1-\bar{I}_2)\bar{\gamma}_1\bar{\gamma}_2}{\bar{I}_1\bar{I}_2} +
gh(\bar{\gamma}_4\chi_2-\bar{\gamma}_5\chi_1) ),
\end{align*}
since $\nabla_{\Pi_i}\Pi_i=1,\; \nabla_{\Pi_i}\Pi_j=0, \; i\neq j, \;  \nabla_{\Pi_i}\Gamma_j=\nabla_{\Gamma_i}\Pi_j=0$
and $\chi=(\chi_1,\chi_2,\chi_3), \; \nabla_{\Gamma_j}
(h_{(\mu,a)})= gh\chi_j, \; \nabla_{\Pi_3} (h_{(\mu,a)}
)= (\Pi_3-l )/\bar{I}_3 ,\; \nabla_{\Pi_k}
(h_{(\mu,a)})= \Pi_k/\bar{I}_k , \; \frac{\partial
\Pi}{\partial \alpha}= \frac{\partial (h_{(\mu,a)})}{\partial
\alpha}=0, \; i, j=1,2,3, \; k=1,2. $

\begin{align*}
 X_{h_{(\mu,a)}}(\Gamma) \cdot \bar{\gamma} & =\{\Gamma,\;
h_{(\mu,a)}\}_{-}|_{\mathcal{O}_{(\mu,a)} \times \mathbb{R}\times
\mathbb{R}^{*}}\cdot \bar{\gamma}\\
& =-\Pi\cdot(\nabla_\Pi\Gamma\times\nabla_\Pi (h_{(\mu,a)})) \cdot
\bar{\gamma}-\Gamma\cdot(\nabla_\Pi\Gamma\times\nabla_\Gamma
(h_{(\mu,a)})-\nabla_\Pi
(h_{(\mu,a)}) \times\nabla_\Gamma\Gamma)\cdot \bar{\gamma}\\
& \;\;\;\; + \{\Gamma,\;
h_{(\mu,a)}\}_{\mathbb{R}}|_{\mathcal{O}_{(\mu,a)} \times
\mathbb{R}\times \mathbb{R}^{*}}\cdot \bar{\gamma}\\
& =\nabla_\Gamma\Gamma\cdot(\Gamma\times\nabla_\Pi
(h_{(\mu,a)}))\cdot \bar{\gamma}+ (\frac{\partial
\Gamma}{\partial \alpha} \frac{\partial (h_{(\mu,a)})}{\partial
l}- \frac{\partial (h_{(\mu,a)})}{\partial
\alpha}\frac{\partial \Gamma}{\partial l})\cdot \bar{\gamma}\\
&=(\Gamma_1,\Gamma_2,\Gamma_3)\times (\frac{ \Pi_1}{\bar{I}_1},
\frac{ \Pi_2}{ \bar{I}_2}, \frac{(\Pi_3-l)}{\bar{I}_3})\cdot \bar{\gamma}\\
&= ( \frac{\bar{I}_2\bar{\gamma}_5\bar{\gamma}_3-\bar{I}_3\bar{\gamma}_6\bar{\gamma}_2-
\bar{I}_2\bar{\gamma}_2\bar{\gamma}_{8}}{\bar{I}_2\bar{I}_3}, \;\;
\frac{\bar{I}_3\bar{\gamma}_6\bar{\gamma}_1-\bar{I}_1\bar{\gamma}_4\bar{\gamma}_3+
\bar{I}_1\bar{\gamma}_1\bar{\gamma}_8}{\bar{I}_3\bar{I}_1}, \;\;
\frac{\bar{I}_1\bar{\gamma}_4\bar{\gamma}_2-
\bar{I}_2\bar{\gamma}_5\bar{\gamma}_1}{\bar{I}_1\bar{I}_2} ),
\end{align*}
since $\nabla_{\Gamma_i}\Gamma_i =1, \; \nabla_{\Gamma_i}\Gamma_j=0, \; i\neq j, \; \nabla_{\Pi_i}\Gamma_j =0, $ and
$\nabla_{\Pi_3} (h_{(\mu,a)})= (\Pi_3-l )/\bar{I}_3 ,\; \nabla_{\Pi_k}
(h_{(\mu,a)})= \Pi_k/\bar{I}_k , \; \frac{\partial
\Gamma_j}{\partial \alpha}= \frac{\partial (h_{(\mu,a)})}{\partial
\alpha}=0, \; i, j= 1,2,3, \; k=1,2.$

\begin{align*}
X_{h_{(\mu,a)}}(\alpha) \cdot \bar{\gamma} & =\{\alpha,\;
h_{(\mu,a)}\}_{-}|_{\mathcal{O}_{(\mu,a)} \times \mathbb{R}\times
\mathbb{R}^{*}}\cdot \bar{\gamma}\\
& =-\Pi\cdot(\nabla_\Pi \alpha \times\nabla_\Pi (h_{(\mu,a)}))\cdot
\bar{\gamma}-\Gamma\cdot(\nabla_\Pi \alpha \times\nabla_\Gamma
(h_{(\mu,a)})-\nabla_\Pi
(h_{(\mu,a)}) \times\nabla_\Gamma \alpha)\cdot \bar{\gamma}\\
& \;\;\;\; + \{\alpha,\;
h_{(\mu,a)}\}_{\mathbb{R}}|_{\mathcal{O}_{(\mu,a)} \times
\mathbb{R}\times \mathbb{R}^{*}}\cdot \bar{\gamma}\\
& = (\frac{\partial \alpha}{\partial \alpha}
\frac{\partial (h_{(\mu,a)})}{\partial l}- \frac{\partial
(h_{(\mu,a)})}{\partial
\alpha}\frac{\partial \alpha}{\partial l})\cdot \bar{\gamma}
 = -\frac{(\bar{\gamma}_3- \bar{\gamma}_8)}{ \bar{I}_3}+
\frac{\bar{\gamma}_8}{J_3},
\end{align*}
since $\nabla_{\Pi_i}\alpha=\nabla_{\Gamma_i}\alpha=0, $ $\frac{\partial \alpha}{\partial
\alpha}= 1, \;\; \frac{\partial (h_{(\mu,a)})}{\partial \alpha}=0, $
and $\frac{\partial (h_{(\mu,a)})}{\partial l}= -(\Pi_3-l )/\bar{I}_3
+\frac{l}{J_3}, \; i= 1,2,3. $

\begin{align*}
 X_{h_{(\mu,a)}}(l) \cdot \bar{\gamma} & =\{l,\;
h_{(\mu,a)}\}_{-}|_{\mathcal{O}_{(\mu,a)} \times \mathbb{R}\times
\mathbb{R}^{*}}\cdot \bar{\gamma}\\
& =-\Pi\cdot(\nabla_\Pi l \times\nabla_\Pi (h_{(\mu,a)})) \cdot
\bar{\gamma}-\Gamma\cdot(\nabla_\Pi l \times\nabla_\Gamma
(h_{(\mu,a)})-\nabla_\Pi
(h_{(\mu,a)}) \times\nabla_\Gamma l)\cdot \bar{\gamma}\\
& \;\;\;\; + \{l,\;
h_{(\mu,a)}\}_{\mathbb{R}}|_{\mathcal{O}_{(\mu,a)} \times
\mathbb{R}\times \mathbb{R}^{*}}\cdot \bar{\gamma}\\
& = (\frac{\partial l}{\partial \alpha} \frac{\partial
(h_{(\mu,a)})}{\partial l}- \frac{\partial (h_{(\mu,a)})}{\partial
\alpha}\frac{\partial l}{\partial l})\cdot \bar{\gamma}=0,
\end{align*}
since $\nabla_{\Pi_i} l=\nabla_{\Gamma_i} l=0,$ and $\frac{\partial l}{\partial \alpha}=
\frac{\partial (h_{(\mu,a)})}{\partial \alpha}=0, \; i=1,2,3. $\\

On the other hand, from the expressions of the dynamical vector field $\tilde{X}$
and Hamiltonian vector field $X_H$, we have that
\begin{align*}
\tilde{X}(\Pi, \Gamma, \alpha, l)^\gamma & =T\pi_{\textmd{SE}(3)\times \mathbb{R}}\cdot \tilde{X}\cdot\gamma(\Pi, \Gamma, \alpha, l)\\
& =T\pi_{\textmd{SE}(3)\times \mathbb{R}}\cdot (X_H+ \textnormal{vlift}(u))\cdot\gamma (\Pi, \Gamma, \alpha, l)\\
&=T\pi_{\textmd{SE}(3)\times \mathbb{R}}\cdot X_H \cdot\gamma (\Pi, \Gamma, \alpha, l) = X_H\cdot\gamma(\Pi, \Gamma, \alpha, l),
\end{align*}
that is,
\begin{align*}
\tilde{X}(\Pi)^\gamma & = X_H(\Pi)\cdot\gamma \\
& = ( \frac{(\bar{I}_2-\bar{I}_3)\gamma_8\gamma_9-
\bar{I}_2\gamma_8\gamma_{14}}{\bar{I}_2\bar{I}_3}+
gh(\gamma_{11}\chi_3-\gamma_{12}\chi_2),\\
& \;\;\;\;\;\;
\frac{(\bar{I}_3-\bar{I}_1)\gamma_9\gamma_7+
\bar{I}_1\gamma_7\gamma_{14}}{\bar{I}_3\bar{I}_1}+gh(\gamma_{12}\chi_1-\gamma_{10}\chi_3), \;\;
\frac{(\bar{I}_1-\bar{I}_2)\gamma_7\gamma_8}{\bar{I}_1\bar{I}_2}
 + gh(\gamma_{10}\chi_2-\gamma_{11}\chi_1) ),
\end{align*}
\begin{align*}
\tilde{X}(\Gamma)^\gamma & = X_H(\Gamma)\cdot\gamma \\
&= ( \frac{\bar{I}_2\gamma_{11}\gamma_9-\bar{I}_3\gamma_{12}\gamma_8-
\bar{I}_2\gamma_8\gamma_{14}}{\bar{I}_2\bar{I}_3}, \;\;
\frac{\bar{I}_3\gamma_{12}\gamma_7-\bar{I}_1\gamma_{10}\gamma_9+
\bar{I}_1\gamma_7\gamma_{14}}{\bar{I}_3\bar{I}_1}, \;\;
\frac{\bar{I}_1\gamma_{10}\gamma_8-\bar{I}_2\gamma_{11}\gamma_7}{\bar{I}_1\bar{I}_2} ),
\end{align*}
\begin{align*}
\tilde{X}(\alpha)^\gamma = X_H(\alpha)\cdot\gamma
 =  -\frac{(\gamma_9- \gamma_{14})}{ \bar{I}_3}
+\frac{\gamma_{14}}{J_3},\;\;\;\;\;\;\;\;
\tilde{X}(l)^\gamma  = X_H(l)\cdot\gamma= 0.
\end{align*}
Since $\gamma$ is closed with respect to
$T\pi_{(\textmd{SE}(3) \times \mathbb{R})}: TT^* (\textmd{SE}(3) \times \mathbb{R})
\rightarrow T(\textmd{SE}(3) \times \mathbb{R}), $
then $\pi_{(\textmd{SE}(3) \times \mathbb{R})}^*(\mathbf{d}\gamma)=0.$ We choose that
$(\gamma_7,\gamma_8,\gamma_9)=\Pi=(\Pi_1,\Pi_2,\Pi_3)=
(\bar{\gamma}_1,\bar{\gamma}_2,\bar{\gamma}_3), \;
(\gamma_{10},\gamma_{11},\gamma_{12})=\Gamma=(\Gamma_1,\Gamma_2,\Gamma_3)=
(\bar{\gamma}_4,\bar{\gamma}_5,\bar{\gamma}_6), $
and $ \gamma_{13}= \alpha= \bar{\gamma}_7,
\;  \gamma_{14}= l = \bar{\gamma}_8. $ Hence
\begin{align*}
& T\bar{\gamma}\cdot \tilde{X}(\Pi)^\gamma
= X_{h_{(\mu,a)}}(\Pi) \cdot \bar{\gamma}, \;\;\;\;\;\;
T\bar{\gamma}\cdot \tilde{X}(\Gamma)^\gamma
= X_{h_{(\mu,a)}}(\Gamma) \cdot \bar{\gamma}, \\
& T\bar{\gamma}\cdot \tilde{X}(\alpha)^\gamma
= X_{h_{(\mu,a)}}(\alpha) \cdot \bar{\gamma}, \;\;\;\;\;\;
T\bar{\gamma}\cdot \tilde{X}(l)^\gamma
= X_{h_{(\mu,a)}}(l) \cdot \bar{\gamma}.
\end{align*}
Thus, the Type I of Hamilton-Jacobi equation for the
$R_p$-reduced controlled rigid spacecraft-rotor system
$(\mathcal{O}_{(\mu,a)}\times \mathbb{R}\times \mathbb{R}^{*},
\omega^{-}_{\mathcal{O}_{(\mu,a)}\times \mathbb{R}\times \mathbb{R}^{*}},
h_{(\mu,a)}, u_{(\mu,a)})$ holds.\\

Next, for $(\mu,a) \in \mathfrak{se}^\ast(3),$
the regular value of $\mathbf{J}_Q$, and a $\textmd{SE}(3)_{(\mu,a)}$-invariant symplectic map
$\varepsilon: T^* (\textmd{SE}(3)\times \mathbb{R}) \rightarrow T^*(\textmd{SE}(3)\times \mathbb{R}),$
assume that $\varepsilon(A,c,\Pi,\Gamma, \alpha, l) =(\varepsilon_1,\cdots,
\varepsilon_{14}),$ and
$\varepsilon(\mathbf{J}^{-1}((\mu,a)))\subset
\mathbf{J}^{-1}((\mu,a)). $ Denote by
$\bar{\varepsilon}=\pi_{(\mu,a)}(\varepsilon):
\mathbf{J}^{-1}((\mu,a)) \rightarrow \mathcal{O}_{(\mu,a)}\times \mathbb{R}\times \mathbb{R}^{*}, $
and $\bar{\varepsilon}=(\bar{\varepsilon}_1,\cdots,\bar{\varepsilon}_{8})
\in \mathcal{O}_{(\mu,a)}\times \mathbb{R}\times \mathbb{R}^{*}
(\subset \mathfrak{se}^\ast(3)\times \mathbb{R}\times \mathbb{R}^{*}), $
and $\lambda= \gamma \cdot \pi_{\textmd{SE}(3)\times \mathbb{R}}: T^* (\textmd{SE}(3)\times \mathbb{R})
\rightarrow T^* (\textmd{SE}(3)\times \mathbb{R}),$ and $\lambda(A,c,\Gamma, \Pi, \alpha, l)
=(\lambda_1,\cdots, \lambda_{14}),$ and
$\bar{\lambda}=\pi_{(\mu,a)}(\lambda): T^* (\textmd{SE}(3)\times \mathbb{R})
\rightarrow \mathcal{O}_{(\mu,a)}\times \mathbb{R}\times \mathbb{R}^{*}, $ and
$\bar{\lambda}=
(\bar{\lambda}_1,\cdots,\bar{\lambda}_{8}) \in \mathcal{O}_{(\mu,a)}\times \mathbb{R}\times \mathbb{R}^{*}
(\subset \mathfrak{se}^\ast(3)\times \mathbb{R}\times \mathbb{R}^{*}). $ We choose that
$(\Pi,\Gamma,\alpha,l)\in
\mathcal{O}_{(\mu,a)}\times \mathbb{R} \times \mathbb{R}^{*}, $ and
$\Pi=(\Pi_1,\Pi_2,\Pi_3)
=(\bar{\varepsilon}_1,\bar{\varepsilon}_2,\bar{\varepsilon}_3),
\; \Gamma= (\Gamma_1,\Gamma_2,\Gamma_3)=
(\bar{\varepsilon}_4,\bar{\varepsilon}_5,\bar{\varepsilon}_6), $
$ \alpha = \bar{\varepsilon}_7, $ and $ l =\bar{\varepsilon}_8 $.
Then $h_{(\mu,a)} \cdot \bar{\varepsilon}:
T^*(\textmd{SE}(3)\times \mathbb{R}) \rightarrow \mathbb{R} $ is given by
\begin{align*}
& h_{(\mu,a)}(\Pi,\Gamma,\alpha,l) \cdot \bar{\varepsilon}=
H(A,c,\Pi,\Gamma, \alpha, l)|_{\mathcal{O}_{(\mu,a)}\times \mathbb{R}\times \mathbb{R}^{*}} \cdot \bar{\varepsilon}\\
& =\frac{1}{2}[\frac{ \bar{\varepsilon}_1^2}{\bar{I}_1}+\frac{ \bar{\varepsilon}_2^2}{\bar{I}_2}
+\frac{(\bar{\varepsilon}_3-\bar{\varepsilon}_8)^2}{\bar{I}_3}+\frac{\bar{\varepsilon}_8^2}{J_3}]
+ gh(\bar{\varepsilon}_4\cdot\chi_1+ \bar{\varepsilon}_5\cdot\chi_2+
\bar{\varepsilon}_6\cdot\chi_3),
\end{align*} and the vector field
\begin{align*}
& X_{h_{(\mu,a)}}(\Pi) \cdot \bar{\varepsilon}
= \{\Pi,h_{(\mu,a)}\}_{-}|_{\mathcal{O}_{(\mu,a)}\times \mathbb{R}\times \mathbb{R}^{*}}\cdot \bar{\varepsilon}\\
&=(\Pi_1,\Pi_2,\Pi_3)\times (\frac{ \Pi_1}{ \bar{I}_1},
\frac{ \Pi_2}{ \bar{I}_2}, \frac{(\Pi_3-l)}{\bar{I}_3})\cdot
\bar{\varepsilon}
+gh(\Gamma_1,\Gamma_2,\Gamma_3)\times (\chi_1,\chi_2,\chi_3)\cdot \bar{\varepsilon}\\
&= ( \frac{(\bar{I}_2-\bar{I}_3)\bar{\varepsilon}_2\bar{\varepsilon}_3-
\bar{I}_2\bar{\varepsilon}_2\bar{\varepsilon}_{8}}{\bar{I}_2\bar{I}_3}+
gh(\bar{\varepsilon}_5\chi_3-\bar{\varepsilon}_6\chi_2), \\
& \;\;\;\;\;\; \frac{(\bar{I}_3-\bar{I}_1)\bar{\varepsilon}_3
\bar{\varepsilon}_1+
\bar{I}_1\bar{\varepsilon}_1\bar{\varepsilon}_8}{\bar{I}_3\bar{I}_1} +
gh(\bar{\varepsilon}_6\chi_1-\bar{\varepsilon}_4\chi_3), \;\;
\frac{(\bar{I}_1-\bar{I}_2)\bar{\varepsilon}_1\bar{\varepsilon}_2}{\bar{I}_1\bar{I}_2} +
gh(\bar{\varepsilon}_4\chi_2-\bar{\varepsilon}_5\chi_1) ),
\end{align*}
\begin{align*}
& X_{h_{(\mu,a)}}(\Gamma) \cdot
\bar{\varepsilon}
= \{\Gamma,h_{(\mu,a)}\}_{-}|_{\mathcal{O}_{(\mu,a)}\times \mathbb{R}\times \mathbb{R}^{*}}\cdot \bar{\varepsilon}
=(\Gamma_1,\Gamma_2,\Gamma_3)\times (\frac{ \Pi_1}{
\bar{I}_1}, \frac{ \Pi_2}{ \bar{I}_2},
\frac{(\Pi_3-l)}{\bar{I}_3})\cdot \bar{\varepsilon}\\
&= ( \frac{\bar{I}_2\bar{\varepsilon}_5\bar{\varepsilon}_3-\bar{I}_3\bar{\varepsilon}_6\bar{\varepsilon}_2-
\bar{I}_2\bar{\varepsilon}_2\bar{\varepsilon}_{8}}{\bar{I}_2\bar{I}_3}, \;\;
\frac{\bar{I}_3\bar{\varepsilon}_6\bar{\varepsilon}_1-\bar{I}_1\bar{\varepsilon}_4\bar{\varepsilon}_3-
\bar{I}_1\bar{\varepsilon}_1\bar{\varepsilon}_8}{\bar{I}_3\bar{I}_1}, \;\;
\frac{\bar{I}_1\bar{\varepsilon}_4\bar{\varepsilon}_2-\bar{I}_2\bar{\varepsilon}_5\bar{\varepsilon}_1}{\bar{I}_1\bar{I}_2} ),
\end{align*}
\begin{align*}
 X_{h_{(\mu,a)}}(\alpha) \cdot \bar{\varepsilon}& =\{\alpha,\;
h_{(\mu,a)}\}_{-}|_{\mathcal{O}_{(\mu,a)} \times \mathbb{R}\times
\mathbb{R}^{*}}\cdot \bar{\varepsilon}\\
& = (\frac{\partial \alpha}{\partial \alpha}
\frac{\partial (h_{(\mu,a)})}{\partial l}- \frac{\partial
(h_{(\mu,a)})}{\partial \alpha }\frac{\partial \alpha}{\partial l})\cdot \bar{\varepsilon}
= -\frac{(\bar{\varepsilon}_3- \bar{\varepsilon}_8)}{ \bar{I}_3}+
\frac{\bar{\varepsilon}_8}{J_3},
\end{align*}
\begin{align*}
& X_{h_{(\mu,a)}}(l) \cdot \bar{\varepsilon}=\{l,\;
h_{(\mu,a)}\}_{-}|_{\mathcal{O}_{(\mu,a)} \times \mathbb{R}\times
\mathbb{R}^{*}}\cdot \bar{\varepsilon}
= (\frac{\partial l}{\partial \alpha} \frac{\partial
(h_{(\mu,a)})}{\partial l}- \frac{\partial (h_{(\mu,a)})}{\partial
\alpha}\frac{\partial l}{\partial l})\cdot \bar{\varepsilon}= 0.
\end{align*}
On the other hand, from the expressions of the dynamical vector field $\tilde{X}$
and Hamiltonian vector field $X_H$, we have that
\begin{align*}
\tilde{X}(\Pi, \Gamma, \alpha, l)^\varepsilon & =T\pi_{\textmd{SE}(3)\times \mathbb{R}}\cdot \tilde{X}\cdot\varepsilon(\Pi, \Gamma, \alpha, l)\\
& =T\pi_{\textmd{SE}(3)\times \mathbb{R}}\cdot (X_H+ \textnormal{vlift}(u))\cdot\varepsilon (\Pi, \Gamma, \alpha, l)\\
& =T\pi_{\textmd{SE}(3)\times \mathbb{R}}\cdot X_H \cdot\varepsilon (\Pi, \Gamma, \alpha, l)= X_H\cdot\varepsilon(\Pi, \Gamma, \alpha, l),
\end{align*}
that is,
\begin{align*}
\tilde{X}(\Pi)^\varepsilon & = X_H(\Pi)\cdot\varepsilon \\
& = ( \frac{(\bar{I}_2-\bar{I}_3)\varepsilon_8\varepsilon_9-
\bar{I}_2\varepsilon_8\varepsilon_{14}}{\bar{I}_2\bar{I}_3}+
gh(\varepsilon_{11}\chi_3-\varepsilon_{12}\chi_2),\\
& \;\;\;\;\;\;
\frac{(\bar{I}_3-\bar{I}_1)\varepsilon_9\varepsilon_7+
\bar{I}_1\varepsilon_7\varepsilon_{14}}{\bar{I}_3\bar{I}_1}
+gh(\varepsilon_{12}\chi_1-\varepsilon_{10}\chi_3), \;\;
\frac{(\bar{I}_1-\bar{I}_2)\varepsilon_7\varepsilon_8}{\bar{I}_1\bar{I}_2}
 + gh(\varepsilon_{10}\chi_2-\varepsilon_{11}\chi_1) ),
\end{align*}
\begin{align*}
\tilde{X}(\Gamma)^\varepsilon & = X_H(\Gamma)\cdot\varepsilon \\
&= ( \frac{\bar{I}_2\varepsilon_{11}\varepsilon_9-\bar{I}_3\varepsilon_{12}\varepsilon_8-
\bar{I}_2\varepsilon_8\varepsilon_{14}}{\bar{I}_2\bar{I}_3}, \;\;
\frac{\bar{I}_3\varepsilon_{12}\varepsilon_7-\bar{I}_1\varepsilon_{10}\varepsilon_9+
\bar{I}_1\varepsilon_7\varepsilon_{14}}{\bar{I}_3\bar{I}_1}, \;\;
\frac{\bar{I}_1\varepsilon_{10}\varepsilon_8-\bar{I}_2\varepsilon_{11}\varepsilon_7}{\bar{I}_1\bar{I}_2} ),
\end{align*}
\begin{align*}
\tilde{X}(\alpha)^\varepsilon  = X_H(\alpha)\cdot\varepsilon
= -\frac{(\varepsilon_9- \varepsilon_{14})}{ \bar{I}_3}
+\frac{\varepsilon_{14}}{J_3}, \;\;\;\;\;\;\;\;
\tilde{X}(l)^\varepsilon  = X_H(l)\cdot\varepsilon =0,
\end{align*}
then we have that
\begin{align*}
T\bar{\gamma}\cdot \tilde{X}(\Pi)^\varepsilon
&= ( \frac{(\bar{I}_2-\bar{I}_3)\bar{\gamma}_2\bar{\gamma}_3-
\bar{I}_2\bar{\gamma}_2\bar{\gamma}_{8}}{\bar{I}_2\bar{I}_3}+
gh(\bar{\gamma}_5\chi_3-\bar{\gamma}_6\chi_2),\\
& \;\;\;\;\;\; \frac{(\bar{I}_3-\bar{I}_1)\bar{\gamma}_3
\bar{\gamma}_1+
\bar{I}_1\bar{\gamma}_1\bar{\gamma}_8}{\bar{I}_3\bar{I}_1} +
gh(\bar{\gamma}_6\chi_1-\bar{\gamma}_4\chi_3), \;\;\;
\frac{(\bar{I}_1-\bar{I}_2)\bar{\gamma}_1\bar{\gamma}_2}{\bar{I}_1\bar{I}_2} +
gh(\bar{\gamma}_4\chi_2-\bar{\gamma}_5\chi_1) ),
\end{align*}
\begin{align*}
T\bar{\gamma}\cdot \tilde{X}(\Gamma)^\varepsilon
&= ( \frac{\bar{I}_2\bar{\gamma}_5\bar{\gamma}_3-\bar{I}_3\bar{\gamma}_6\bar{\gamma}_2-
\bar{I}_2\bar{\gamma}_2\bar{\gamma}_{8}}{\bar{I}_2\bar{I}_3}, \;\;
\frac{\bar{I}_3\bar{\gamma}_6\bar{\gamma}_1-\bar{I}_1\bar{\gamma}_4\bar{\gamma}_3+
\bar{I}_1\bar{\gamma}_1\bar{\gamma}_8}{\bar{I}_3\bar{I}_1}, \;\;
\frac{\bar{I}_1\bar{\gamma}_4\bar{\gamma}_2-
\bar{I}_2\bar{\gamma}_5\bar{\gamma}_1}{\bar{I}_1\bar{I}_2} ),
\end{align*}
\begin{align*}
T\bar{\gamma}\cdot \tilde{X}(\alpha)^\varepsilon
 = -\frac{(\bar{\gamma}_3- \bar{\gamma}_8)}{ \bar{I}_3}+
\frac{\bar{\gamma}_8}{J_3}, \;\;\;\;\;\;\;\;
T\bar{\gamma}\cdot \tilde{X}(l)^\varepsilon= 0.
\end{align*}
Note that $$T\bar{\lambda}\cdot \tilde{X} \cdot \varepsilon=T\pi_{(\mu,a)}\cdot T\lambda \cdot (X_H+ \textnormal{vlift}(u))\cdot\varepsilon
=T\pi_{(\mu,a)}\cdot T\gamma \cdot T\pi_{\textmd{SE}(3)\times \mathbb{R}}\cdot (X_H+ \textnormal{vlift}(u))\cdot\varepsilon
=T\bar{\lambda}\cdot X_H \cdot \varepsilon,$$ that is,
\begin{align*}
T\bar{\lambda}\cdot \tilde{X}(\Pi) \cdot \varepsilon & =T\bar{\lambda}\cdot X_H(\Pi) \cdot \varepsilon\\
&= ( \frac{(\bar{I}_2-\bar{I}_3)\bar{\lambda}_2\bar{\lambda}_3-
\bar{I}_2\bar{\lambda}_2\bar{\lambda}_{8}}{\bar{I}_2\bar{I}_3}+
gh(\bar{\lambda}_5\chi_3-\bar{\lambda}_6\chi_2),\\
& \;\;\;\;\;\; \frac{(\bar{I}_3-\bar{I}_1)\bar{\lambda}_3
\bar{\lambda}_1+
\bar{I}_1\bar{\lambda}_1\bar{\lambda}_8}{\bar{I}_3\bar{I}_1} +
gh(\bar{\lambda}_6\chi_1-\bar{\lambda}_4\chi_3), \;\;
\frac{(\bar{I}_1-\bar{I}_2)\bar{\lambda}_1\bar{\lambda}_2}{\bar{I}_1\bar{I}_2} +
gh(\bar{\lambda}_4\chi_2-\bar{\lambda}_5\chi_1) ),
\end{align*}
\begin{align*}
T\bar{\lambda}\cdot \tilde{X}(\Gamma) \cdot \varepsilon & =T\bar{\lambda}\cdot X_H(\Gamma) \cdot \varepsilon\\
&= ( \frac{\bar{I}_2\bar{\lambda}_5\bar{\lambda}_3-\bar{I}_3\bar{\lambda}_6\bar{\lambda}_2-
\bar{I}_2\bar{\lambda}_2\bar{\lambda}_{8}}{\bar{I}_2\bar{I}_3}, \;\;
\frac{\bar{I}_3\bar{\lambda}_6\bar{\lambda}_1-\bar{I}_1\bar{\lambda}_4\bar{\lambda}_3+
\bar{I}_1\bar{\lambda}_1\bar{\lambda}_8}{\bar{I}_3\bar{I}_1}, \;\;
\frac{\bar{I}_1\bar{\lambda}_4\bar{\lambda}_2-
\bar{I}_2\bar{\lambda}_5\bar{\lambda}_1}{\bar{I}_1\bar{I}_2} ),
\end{align*}
\begin{align*}
T\bar{\lambda}\cdot \tilde{X}(\alpha) \cdot \varepsilon=T\bar{\lambda}\cdot X_H(\alpha) \cdot \varepsilon
= -\frac{(\bar{\lambda}_3- \bar{\lambda}_8)}{ \bar{I}_3}+
\frac{\bar{\lambda}_8}{J_3}, \;\;\;\;\;\;\;\;
T\bar{\lambda}\cdot \tilde{X}(l) \cdot \varepsilon=T\bar{\lambda}\cdot X_H(l) \cdot \varepsilon= 0.
\end{align*}
Thus, when we choose that
$(\Pi,\Gamma,\alpha,l)\in
\mathcal{O}_{(\mu,a)}\times \mathbb{R}\times \mathbb{R}^{*}, $ and
 $(\varepsilon_7,\varepsilon_8,\varepsilon_9)=\Pi=(\Pi_1,\Pi_2,\Pi_3)
 =(\bar{\gamma}_1,\bar{\gamma}_2,\bar{\gamma}_3)=
(\bar{\varepsilon}_1,\bar{\varepsilon}_2,\bar{\varepsilon}_3)=
(\bar{\lambda}_1,\bar{\lambda}_2,\bar{\lambda}_3), $ and
$(\varepsilon_{10},\varepsilon_{11},\varepsilon_{12})=\Gamma
=(\Gamma_1,\Gamma_2,\Gamma_3)=(\bar{\gamma}_4,\bar{\gamma}_5,\bar{\gamma}_6)=
(\bar{\varepsilon}_4,\bar{\varepsilon}_5,\bar{\varepsilon}_6)=
(\bar{\lambda}_4,\bar{\lambda}_5,\bar{\lambda}_6), $
and $\varepsilon_{13}=\alpha=\bar{\gamma}_7 =\bar{\varepsilon}_7= \bar{\lambda}_7,
\; \varepsilon_{14}=l=\bar{\gamma}_8= \bar{\varepsilon}_8= \bar{\lambda}_8, $ we must have that
\begin{align*}
& T\bar{\gamma}\cdot \tilde{X}(\Pi)^\varepsilon =X_{h_{(\mu,a)}}(\Pi) \cdot \bar{\varepsilon}
=T\bar{\lambda}\cdot \tilde{X}(\Pi) \cdot \varepsilon, \\
& T\bar{\gamma}\cdot \tilde{X}(\Gamma)^\varepsilon=X_{h_{(\mu,a)}}(\Gamma) \cdot \bar{\varepsilon}
=T\bar{\lambda}\cdot \tilde{X}(\Gamma) \cdot \varepsilon,\\
& T\bar{\gamma}\cdot \tilde{X}(\alpha)^\varepsilon =X_{h_{(\mu,a)}}(\alpha) \cdot \bar{\varepsilon}
=T\bar{\lambda}\cdot \tilde{X}(\alpha) \cdot \varepsilon, \\
& T\bar{\gamma}\cdot \tilde{X}(l)^\varepsilon =X_{h_{(\mu,a)}}(l) \cdot \bar{\varepsilon}
=T\bar{\lambda}\cdot \tilde{X}(l) \cdot \varepsilon.
\end{align*}
Since the map $\varepsilon: T^* (\textmd{SE}(3)\times \mathbb{R})
\rightarrow T^* (\textmd{SE}(3)\times \mathbb{R})$ is symplectic, then
$T\bar{\varepsilon}\cdot X_{h_{(\mu,a)} \cdot \bar{\varepsilon}}
=X_{h_{(\mu,a)}} \cdot \bar{\varepsilon}. $
Thus, in this case, we must have that
$\varepsilon$ and $\bar{\varepsilon} $ are the solution of the Type II of
Hamilton-Jacobi equation
$T\bar{\gamma}\cdot \tilde{X}^\varepsilon= X_{h_{(\mu,a)}}\cdot \bar{\varepsilon}, $
for the $R_p$-reduced controlled rigid spacecraft-rotor system
$(\mathcal{O}_{(\mu,a)}\times \mathbb{R}\times \mathbb{R}^{*},
\omega^{-}_{\mathcal{O}_{(\mu,a)}\times \mathbb{R}\times \mathbb{R}^{*}},
h_{(\mu,a)}, u_{(\mu,a)})$, if and only if they satisfy
the equation $T\bar{\varepsilon}\cdot(X_{h_{(\mu,a)} \cdot \bar{\varepsilon}})
= T\bar{\lambda}\cdot \tilde{X}\cdot\varepsilon. $\\

To sum up the above discussion, we have the following Theorem 4.5.
For convenience, the maps involved in
the following theorem are shown in Diagram-3.

\begin{center}
\hskip 0cm \xymatrix{ \mathbf{J}_Q^{-1}(\mu,a) \ar[r]^{i_{(\mu,a)}} & T^* Q
\ar[d]_{X_{H\cdot \varepsilon}} \ar[dr]^{\tilde{X}^\varepsilon} \ar[r]^{\pi_Q}
& Q \ar[d]^{\tilde{X}^\gamma} \ar[r]^{\gamma}
& T^*Q \ar[d]_{\tilde{X}} \ar[dr]_{X_{h_{(\mu,a)} \cdot\bar{\varepsilon}}} \ar[r]^{\pi_{(\mu,a)}}
& \;\;\; \mathcal{O}_{(\mu,a)}\times \mathbb{R}\times \mathbb{R}^{*} \ar[d]^{X_{h_{(\mu,a)}}} \\
& T(T^*Q)  & TQ \ar[l]^{T\gamma}
& T(T^*Q) \ar[l]^{T\pi_Q} \ar[r]_{T\pi_{(\mu,a)}}
& \;\;\; T(\mathcal{O}_{(\mu,a)}\times \mathbb{R}\times \mathbb{R}^{*})}
\end{center}
$$\mbox{Diagram-3}$$

\begin{theo}
In the case of non-coincident centers of buoyancy and gravity,
if the 5-tuple $(T^\ast Q, \textmd{SE}(3),\omega_Q,H,u), $ where $Q=
\textmd{SE}(3)\times \mathbb{R}, $ is a regular point reducible
rigid spacecraft-rotor system with the control torque $u$ acting on the rotor,
then for a point $(\mu,a) \in \mathfrak{se}^\ast(3)$, the regular
value of the momentum map $\mathbf{J}_Q: \textmd{SE}(3)\times
\mathfrak{se}^\ast(3) \times \mathbb{R}\times \mathbb{R}^{*} \to
\mathfrak{se}^\ast(3)$, the $R_p$-reduced controlled rigid spacecraft-rotor system is the 4-tuple
$(\mathcal{O}_{(\mu,a)} \times \mathbb{R}\times
\mathbb{R}^{*},\tilde{\omega}_{\mathcal{O}_{(\mu,a)} \times \mathbb{R}
\times \mathbb{R}^{*}}^{-},h_{(\mu,a)},u_{(\mu,a)}). $ Assume that
$\gamma: \textmd{SE}(3)\times \mathbb{R}\rightarrow
T^*(\textmd{SE}(3)\times \mathbb{R})$ is an one-form on
$\textmd{SE}(3)\times \mathbb{R}$,
and $\lambda=\gamma \cdot \pi_{(\textmd{SE}(3)\times \mathbb{R})}:
T^* (\textmd{SE}(3)\times \mathbb{R} )\rightarrow T^* (\textmd{SE}(3)\times \mathbb{R}), $ and $\varepsilon:
T^* (\textmd{SE}(3) \times \mathbb{R})\rightarrow T^* (\textmd{SE}(3)\times \mathbb{R}) $ is a
$\textmd{SE}(3)_{(\mu,a)}$-invariant symplectic map.
Denote
$\tilde{X}^\gamma = T\pi_{(\textmd{SE}(3)\times \mathbb{R})}\cdot \tilde{X}\cdot \gamma$, and
$\tilde{X}^\varepsilon = T\pi_{(\textmd{SE}(3)\times \mathbb{R})}\cdot \tilde{X}\cdot \varepsilon$,
where $\tilde{X}=X_{(T^\ast Q,\textmd{SE}(3),\omega_Q,H,u)}$ is the dynamical vector field of
the controlled rigid spacecraft-rotor system $(T^\ast Q,\textmd{SE}(3),\omega_Q,H,u)$.
Moreover, assume that $\textmd{Im}(\gamma)\subset \mathbf{J}_Q^{-1}(\mu,a), $ and it is
$\textmd{SE}(3)_{(\mu,a)}$-invariant,
and $\varepsilon(\mathbf{J}_Q^{-1}(\mu,a))\subset \mathbf{J}_Q^{-1}(\mu,a). $
Denote $\bar{\gamma}=\pi_{(\mu,a)}(\gamma):
\textmd{SE}(3)\times \mathbb{R} \rightarrow \mathcal{O}_{(\mu,a)}\times \mathbb{R}\times \mathbb{R}^{*}, $ and
$\bar{\lambda}=\pi_{(\mu,a)}(\lambda): T^* (\textmd{SE}(3)\times \mathbb{R}) \rightarrow
\mathcal{O}_{(\mu,a)}\times \mathbb{R}\times \mathbb{R}^{*}, $ and
$\bar{\varepsilon}=\pi_{(\mu,a)}(\varepsilon): \mathbf{J}_Q^{-1}(\mu,a)\rightarrow
\mathcal{O}_{(\mu,a)}\times \mathbb{R}\times \mathbb{R}^{*}. $
Then the following two assertions hold:\\
\noindent $(\mathbf{i})$
If the one-form $\gamma: \textmd{SE}(3)\times \mathbb{R} \rightarrow T^*(\textmd{SE}(3)\times \mathbb{R}) $ is closed with respect to
$T\pi_{(\textmd{SE}(3)\times \mathbb{R})}: TT^* (\textmd{SE}(3)\times \mathbb{R}) \rightarrow T(\textmd{SE}(3)\times \mathbb{R}), $
then $\bar{\gamma}$ is a solution of the Type I of Hamilton-Jacobi equation
$T\bar{\gamma}\cdot \tilde{X}^\gamma= X_{h_{(\mu,a)}}\cdot \bar{\gamma}; $\\
\noindent $(\mathbf{ii})$
The $\varepsilon$ and $\bar{\varepsilon} $ satisfy the Type II of Hamilton-Jacobi equation
$T\bar{\gamma}\cdot \tilde{X}^\varepsilon= X_{h_{(\mu,a)}}\cdot \bar{\varepsilon}, $
if and only if they satisfy
the equation $T\bar{\varepsilon}\cdot(X_{h_{(\mu,a)} \cdot \bar{\varepsilon}})
= T\bar{\lambda}\cdot \tilde{X}\cdot\varepsilon. $ \hskip 0.3cm $\blacksquare$
\end{theo}

\begin{rema}
When the rigid spacecraft does not carry any internal rotor, in this
case the configuration space is $Q=G=\textmd{SE}(3), $
the motion of
rigid spacecraft is just the rotation motion with drift of a rigid
body, the above $R_p$-reduced controlled rigid spacecraft-rotor system
is just the Marsden-Weinstein reduced heavy top system, that is,
3-tuple $(\mathcal{O}_{(\mu,a)},
\omega_{\mathcal{O}_{(\mu,a)}},h_{\mathcal{O}_{(\mu,a)}})$, where
$\mathcal{O}_{(\mu,a)} \subset \mathfrak{se}^\ast(3)$ is the
co-adjoint orbit, $\omega_{\mathcal{O}_{(\mu,a)}}$ is orbit
symplectic form on $\mathcal{O}_{(\mu,a)}$, which is induced by the
heavy top Lie-Poisson bracket on $\mathfrak{se}^\ast(3)$,
$h_{\mathcal{O}_{(\mu,a)}}(\Pi,\Gamma)\cdot \pi_{\mathcal{O}_{(\mu,a)}}
=H(A,c,\Pi,\Gamma)|_{\mathcal{O}_{(\mu,a)}}$. From the above Theorem 4.5
we can obtain the Proposition 5.5 in Wang \cite{wa17}, that is, we give the two
types of Lie-Poisson Hamilton-Jacobi equation for the Marsden-Weinstein
reduced heavy top system $(\mathcal{O}_{(\mu,a)},
\omega_{\mathcal{O}_{(\mu,a)}},h_{\mathcal{O}_{(\mu,a)}})$.
See Marsden and Ratiu
\cite{mara99}, Ge and Marsden \cite{gema88}, and Wang \cite{wa17}.
\end{rema}

It is worthy of noting that, for the controlled rigid
spacecraft-rotor system $(T^\ast Q,\textmd{SE}(3),\omega_Q,\\ H,u)$
with the $R_p$-reduced controlled rigid spacecraft-rotor system
$(\mathcal{O}_{(\mu,a)} \times \mathbb{R} \times
\mathbb{R}^{*},\tilde{\omega}_{\mathcal{O}_{(\mu,a)} \times \mathbb{R}
\times \mathbb{R}^{*}}^{-}, \\
h_{(\mu,a)}, u_{(\mu,a)}) $, we know that the Hamiltonian vector fields
$X_{H}$ and $X_{h_{(\mu,a)}}$ for the corresponding
Hamiltonian system $(T^*Q,\textmd{SE}(3),\omega_Q, H)$
and its $R_p$-reduced system $(\mathcal{O}_{(\mu,a)} \times \mathbb{R} \times
\mathbb{R}^{*},\tilde{\omega}_{\mathcal{O}_{(\mu,a)} \times \mathbb{R}
\times \mathbb{R}^{*}}^{-},\\ h_{(\mu,a)} )$, are $\pi_{(\mu,a)}$-related, that is,
$X_{h_{(\mu,a)}}\cdot \pi_{(\mu,a)}=T\pi_{(\mu,a)}\cdot X_{H}\cdot i_{(\mu,a)}.$ By using
the similar way in proof Theorem 4.4, then we can
prove the following Theorem 4.7, which states the relationship
between the solutions of Type II of Hamilton-Jacobi equations and the
regular point reduction.

\begin{theo}
In the case of non-coincident centers of buoyancy and gravity, for the controlled rigid
spacecraft-rotor system $(T^\ast Q,\textmd{SE}(3),\omega_Q,H,u)$
with the $R_p$-reduced controlled rigid spacecraft-rotor system
$(\mathcal{O}_{(\mu,a)} \times \mathbb{R} \times
\mathbb{R}^{*},\tilde{\omega}_{\mathcal{O}_{(\mu,a)} \times \mathbb{R}
\times \mathbb{R}^{*}}^{-}, h_{(\mu,a)}, u_{(\mu,a)}) $,
assume that $\gamma:
\textmd{SE}(3)\times \mathbb{R} \rightarrow
T^*(\textmd{SE}(3)\times \mathbb{R})$ is an one-form on
$\textmd{SE}(3)\times \mathbb{R}$, and $\varepsilon:
T^* (\textmd{SE}(3) \times \mathbb{R})\rightarrow
T^* (\textmd{SE}(3)\times \mathbb{R}) $ is a
$\textmd{SE}(3)_{(\mu,a)}$-invariant symplectic map,
$\bar{\varepsilon}=\pi_{(\mu,a)}(\varepsilon): \mathbf{J}_Q^{-1}(\mu, a )\rightarrow
\mathcal{O}_{(\mu,a)}\times \mathbb{R}\times \mathbb{R}^{*}. $
Under the hypotheses and notations of Theorem 4.5, then we have that
$\varepsilon$ is a solution of the Type II of Hamilton-Jacobi equation
$T\gamma\cdot \tilde{X}^\varepsilon= X_H\cdot \varepsilon, $ for the
regular point reducible controlled rigid
spacecraft-rotor system $(T^\ast Q,\textmd{SE}(3),\omega_Q,H,u), $ if and only if
$\varepsilon$ and $\bar{\varepsilon} $ satisfy the Type II of Hamilton-Jacobi equation
$T\bar{\gamma}\cdot \tilde{X}^\varepsilon= X_{h_{(\mu,a)}}\cdot \bar{\varepsilon}, $ for the
$R_p$-reduced controlled rigid spacecraft-rotor system
$(\mathcal{O}_{(\mu,a)} \times \mathbb{R} \times
\mathbb{R}^{*},\tilde{\omega}_{\mathcal{O}_{(\mu,a)} \times \mathbb{R}
\times \mathbb{R}^{*}}^{-}, h_{(\mu,a)}, u_{(\mu,a)}) $.
\end{theo}

The theory of controlled mechanical system is a very important
subject, following the theoretical and applied development of
geometric mechanics, a lot of important problems about this subject
are being explored and studied. In this paper, we reveal the deeply internal
relationships of the geometrical structures of phase spaces, the dynamical
vector fields and controls of the controlled rigid spacecraft-rotor system.
It is worthy of noting that, in the cases of coincident and
non-coincident centers of buoyancy and gravity,
the motions of the controlled rigid spacecraft-rotor system are different, and
the configuration spaces, the Hamiltonian functions, the actions of Lie group,
the $R_p$-reduced symplectic forms and the $R_p$-reduced systems of
the controlled rigid spacecraft-rotor system are also different. But,
the two types of Hamilton-Jacobi equations given by calculation in detail
are same, that is, the internal rules are same.
It is the key thought of the researches of geometrical mechanics
of the professor Jerrold E. Marsden to explore and study the deeply internal
relationship between the geometrical structure of phase space and the dynamical
vector field of a mechanical system. It is also our goal of pursuing and inheriting.\\

\end{document}